\documentclass[11pt]{amsart}
\usepackage{amsmath,amsfonts,amsthm,amscd,amssymb,mathrsfs,amssymb,bm}
\usepackage{mathrsfs}
\usepackage{graphicx}
\newtheorem{theorem}{Theorem}

\newtheorem{proposition}[theorem]{Proposition}
\newtheorem{lemma}[theorem]{Lemma}

\newtheorem{remark}[theorem]{Remark}
\newcommand{\aaa}{\alpha}

\newcommand{\CCC}{\Gamma}

\newcommand{\eee}{\epsilon}
\newcommand{\CP}{\mathbb{CP}}
\newcommand{\CC}{\mathbb{C}}

\newcommand{\ZZ}{\mathbb{Z}}

\newcommand{\FF}{\mathbb{F}}

\newcommand{\Sing}{{\rm{Sing}\,}}

\newcommand{\ol}{\overline}
\newcommand{\lra}{\longrightarrow}
\newcommand{\lras}{\,\longrightarrow\,}

\newcommand{\set}{\,|\,}
\newcommand{\proofend}{\hfill$\square$}

\newcommand{\inv}{^{-1}}

\newcommand{\Bs}{{\rm{Bs}}}
\newcommand{\Pic}{{\rm{Pic}}}

\newcommand{\ms}{\mathscr}

\numberwithin{equation}{section}
\numberwithin{theorem}{section}

\begin{document}
\bibliographystyle{alpha} 
\title[]
{Geometry of some twistor spaces
of algebraic dimension one}
\author{Nobuhiro Honda}
\address{Department of Mathematics, Tokyo Institute
of Technology}
\email{honda@math.titech.ac.jp}

\thanks{The author has been partially supported by JSPS KAKENHI Grant Number 24540061.
\\
{\it{Mathematics Subject Classification}} (2010) 53A30}
\begin{abstract}
It is shown that there exists a twistor space on the
$n$-fold connected sum of complex projective planes $n\CP^2$,
 whose algebraic dimension is one and whose general fiber of the algebraic reduction is birational to an elliptic ruled surface or a K3 surface.
The former kind of twistor spaces are constructed over $n\CP^2$ for any $n\ge 5$, while the latter kind of example is constructed over $5\CP^2$.
Both of these seem to be the first such example on $n\CP^2$.
The algebraic reduction in these examples 
is induced by the anti-canonical system of the twistor spaces.
It is also shown that the former kind of twistor spaces contain
a pair of non-normal Hopf surfaces.
\end{abstract}
\maketitle

\section{Introduction}

Let $X$ be a compact complex manifold.
We denoted by $\dim X$
for the complex dimension of $X$.
A basic invariant for  $X$ is 
{\emph {the algebraic dimension}},
which is usually denoted by $a(X)$.
This is defined as the transcendental degree over $\CC$ of the field of 
meromorphic functions on $X$,
and thus roughly measures how many meromorphic functions
exist on $X$.
For any compact complex manifold $X$ we have $0\le a(X)\le \dim X$.
We have $a(X)=0$ iff $X$ has only constant meromorphic functions,
and $a(X) = \dim X$ iff $X$ is bimeromorphically equivalent
to a projective algebraic manifold.
For any $X$, there exists a projective algebraic manifold $Y$
and a surjective meromorphic map $f:X\to Y$ which induces
an isomorphism for the meromorphic function fields of $Y$ and $X$,
so that $\dim Y = a(X)$;
the meromorphic map $f:X\to Y$ is called {\em the algebraic reduction}
of $X$, and is known to be unique under a bimeromorphic equivalence.
Moreover, fibers of $f$ are necessarily connected.

If $a(X) = \dim X -1$, a general fiber of the algebraic reduction
is always an elliptic curve; this is a consequence of 
the fact that the degree of the canonical bundle of 
an algebraic curve is zero only when it is an elliptic curve.
In the case $a(X) = \dim X -2$ and $\dim X>2$,
a possible list of a general fiber of the algebraic reduction
is obtained in \cite[p.\,146]{U75};
basically the surfaces in the list are surfaces with non-positive 
Kodaira dimension.
When $X$ is 3-dimensional and belongs to the class $\ms C$,
structure of a possible general fiber of algebraic reduction is determined by A.\,Fujiki
\cite[p.\,236, Theorem]{F83}.
But when $X\not\in\ms C$, it is not easy to construct
a (non-trivial) example of $X$ with $a(X)=1$ which has a surface
in the list
as a general fiber of the algebraic reduction.

As is noticed by F.\,Campana, C.\,LeBrun, Y.\,S.\,Poon, M.\,Ville 
and others, the so called {\em the twistor spaces}
associated to self-dual metrics on real 4-manifolds provide  examples of 
compact complex threefold $Z$ which satisfies $Z\not\in\ms C$ and $a(Z)=1$.
Some  interesting examples are the twistor spaces of 
a Ricci-flat K\"ahler metric on a complex torus and a K3 surface with the complex orientation reversed;
the algebraic reduction of these twistor spaces is a natural holomorphic projection 
to $\CP^1$, which is differentiably 
a fiber bundle over $\CP^1$ whose fibers are a complex torus
or a K3 surface respectively.
Also, a Hopf surface has a conformally flat structure 
whose twistor space enjoys a similar property.
Further, in the case of Hopf surfaces, 
Fujiki \cite{F00-2} showed that if $a(Z)=1$, a general fiber of its algebraic reduction
is either a Hopf surface or a ruled surface over an elliptic curve.
(We will call the latter as an elliptic ruled surface in the sequel.)
These twistor spaces are 
quotient of $\CP^3$ minus two skew lines by a free linear $\ZZ$-action.

As examples of a different flavor, by Donaldson-Friedman \cite{DF89}, 
LeBrun-Poon \cite{LB92,LP92}, Poon \cite{P92} and others,
if $n\ge 4$, the connected sum $n\CP^2$ of 
$n$ copies of the complex projective planes 
admits a self-dual metric whose twistor space $Z$ satisfies $a(Z)=1$.
For these examples, the algebraic reduction is induced by the natural
 square root $K_Z^{-1/2}$
of the anti-canonical line bundle of $Z$, and from this,
it may be readily seen, with a help of a useful result by 
Pedersen-Poon \cite{PP94}, that a general fiber of the algebraic 
reduction  is a rational surface.
Thus rational surfaces actually occur as a general
fiber of the algebraic reduction of a twistor 
space on $n\CP^2$ if $n\ge 4$.
Also, it is not difficult to see that, 
for any twistor space $Z$ on $4\CP^2$,
other complex surfaces cannot be a general fiber
of the algebraic reduction for $Z$ with $a(Z)=1$.
However, to the best of the author's knowledge,
no example is known so far of
a twistor space of algebraic dimension one 
on $n\CP^2$ having any other surfaces in the list as a general fiber
of the algebraic reduction.

In this article, we show that
{\em there exists a
twistor space $Z$ on $n\CP^2$
satisfying $a(Z) = 1$ whose general fiber of the algebraic
reduction is birational to (i) a ruled surface over an elliptic curve, or (ii) a K3 surface.}
The former type of twistor spaces are constructed over $n\CP^2$ 
for arbitrary $n\ge 5$, while the latter is 
constructed only over $5\CP^2$. 
The algebraic reduction for both of these twistor spaces
is induced by the anti-canonical system of the twistor space, and it is a pencil, having base curves.

We briefly explain characteristic features of these twistor spaces.
The twistor spaces with the pencil of  elliptic ruled
surfaces have a particular pair of  rational curves which are 
base curves of the anti-canonical pencil.
Moreover these are
 {\em double curves} of the elliptic ruled surfaces.
In particular, the surfaces are non-normal.
These surfaces are obtained from their normalization
by identifying each of two disjoint sections of the ruling by an involution
on each.
(This is why the two double curves are not elliptic but rational.)
These twistor spaces admit a $\CC^*$-action, and each member
of the pencil is $\CC^*$-invariant.
There exists a  reducible member of the pencil,
whose irreducible components consist of 
two surfaces which are birational to {\em a Hopf surface}.
The intersection of the two components is a non-singular elliptic curve,
and it is a single orbit of the
$\CC^*$-action.
%
On the other hand, the twistor spaces with the pencil of K3 surfaces do not admit a $\CC^*$-action, and the base locus of the pencil consists
of a pair of chains of curves formed by four rational curves. 
From the second Betti number, these twistor spaces
do not seem to have a direct generalization to the case $n>5$.

At a technical level,
both of these twistor spaces are characterized by the property that
they contain a rational surface of special kinds
whose twice is an anti-canonical divisor of the 
twistor space.
The most difficult part for analyzing the structure 
is to show that under this condition
on the presence of the rational surface,
the anti-canonical system is necessarily a pencil.
This is shown by using a variation formula for the 
Euler characteristic for a line bundle under a blow-up
whose center is a curve.
We apply the formula for a particular line bundle and the rational curves, which lead to the conclusion.

\section{Construction of some rational surfaces}
In this section, we first construct a rational elliptic surface
$S_0$ admitting a non-trivial $\CC^*$-action,
which satisfies $K^2 = 0$.
Next we blowup this surface  at fixed points
of the $\CC^*$-action,
and obtain a rational surface $S$ with $\CC^*$-action which satisfies 
$K^2 = 8-2n$, where $n>4$.
Then we investigate the linear system $|mK_S\inv|$ for any $m>0$.

\begin{figure}
\includegraphics{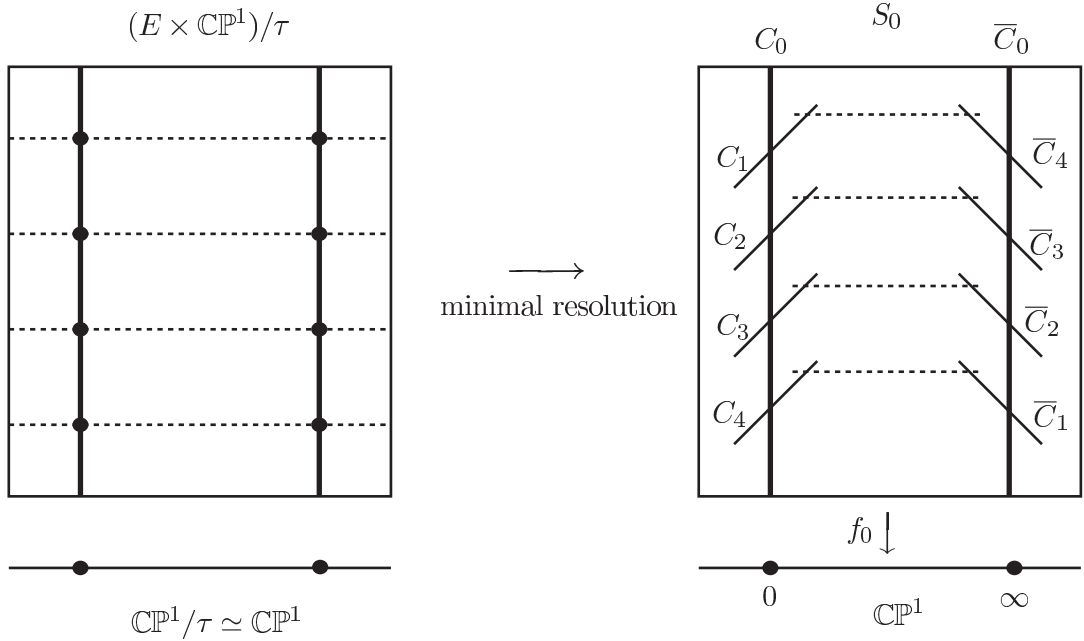}
\caption{Construction of the surface $S_0$ and the singular fibers of the elliptic fibration 
$f_0:S_0\to\CP^1$}
\label{fig:quot}
\end{figure}

By {\em a real structure} on a complex manifold, we mean
an anti-holomorphic involution on it.
Throughout this section $E$ denotes a smooth elliptic curve equipped with:
\begin{itemize}
\item a real structure $\sigma$ without fixed point,
\item a holomorphic involution $\tau$ which has 4 fixed points, and which commutes with $\sigma$.
\end{itemize}
Also on the complex projective line $\CP^1$,  we put:
\begin{itemize}
\item a real structure without fixed point, which is also 
denoted by $\sigma$,
\item a holomorphic involution $\tau$ which has 2 fixed points, and which commutes with $\sigma$.
\end{itemize}
From these, on the product surface $E\times \CP^1$ we obtain:
\begin{itemize}
\item the product real structure, again denoted by $\sigma$,
\item the product holomorphic involution, again denoted by $\tau$,
 commuting with $\sigma$.
\end{itemize}
On the product $E\times\CP^1$, the holomorphic involution $\tau$ has exactly 8 fixed points.
So the quotient complex surface $(E\times\CP^1)/\tau$ has
8 ordinary double points.
This quotient surface  also has a real structure induced by that on $E\times\CP^1$.
The projection $E\times\CP^1\to\CP^1$ induces a holomorphic map
$(E\times\CP^1)/\tau\to\CP^1/\tau$, which defines 
a structure of elliptic surface over $\CP^1/\tau\simeq\CP^1$ on the quotient surface.
Let 
\begin{align}\label{mr0}
S_0 \lras (E\times\CP^1)/\tau
\end{align}
be the minimal resolution of the 8 double points.
The real structure on  $(E\times\CP^1)/\tau$ naturally
lifts to $S_0$, and we still denote it by $\sigma$.
This has no real point.
We note that the surface $S_0$ is uniquely determined by the isomorphism
class of the initial elliptic curve $E$.
Taking the composition with the above projection
$(E\times\CP^1)/\tau\to\CP^1/\tau\simeq\CP^1$,
we obtain an elliptic fibration
\begin{align}\label{ell0}
f_0:S_0\lras \CP^1.
\end{align} 
This has non-reduced fibers 
over the images of the two fixed points of $\tau$
on $\CP^1$.
These fibers are clearly of type $I_0^*$ in Kodaira's notation for singular fibers
of elliptic fibrations.
Denoting $0$ and $\infty$ for the images of the two fixed points,
we write the singular fibers as
\begin{align}\label{I_0*}
f_0\inv(0)= 2C_0 + \sum_{1\le i\le 4} C_i
\quad{\text{and}}\quad
f_0\inv(\infty)=2\ol C_0 + \sum_{1\le i\le 4} \ol C_i,
\end{align}
where $C_0$ and $\ol C_0=\sigma(C_0)$ are the unique double component
of each fiber respectively.
Of course, all of these components are $(-2)$-curves.
(See Figure \ref{fig:quot} for these constructions.)

The following properties on the surface $S_0$ are immediate to see from the construction:

\begin{itemize}
\item The surface $S_0$ has an effective $\CC^*$-action  induced from 
the $\CC^*$-action on $E\times\CP^1$ which is the product of
the trivial action on $E$ and the standard $\CC^*$-action on $\CP^1$.
\item
This $\CC^*$-action on $S_0$ is compatible with the above elliptic fibration
$f_0:S_0\to\CP^1$, and induces a non-trivial but non-effective
$\CC^*$-action on the base space $\CP^1$.
\item All regular fibers of $f_0$ are isomorphic to $E$, and are mutually identified by the $\CC^*$-action.
The two singular fibers of $f_0$ are $\CC^*$-invariant, and
among their components only the
double components $C_0$ and $\ol C_0$ are point-wise fixed
by the $\CC^*$-action.
\end{itemize}

Moreover we have the following properties on 
the structure of the surface $S_0$:

\begin{proposition}\label{prop:rs1}
Let $S_0$ be as above.
Then 
$K^2_{S_0}=0$, and the anti-canonical system 
of the surface $S_0$ is 
a pencil without a fixed point, whose associated morphism
may be identified with the 
elliptic fibration $f_0:S_0\to\CP^1$ in \eqref{ell0}.
Moreover, $S_0$ is a rational surface.
%
%
\end{proposition}

As this is not difficult to see, we omit a proof.
We mention that this surface $S_0$ is the same as
the surface $S'_0$ in \cite[p.\,146]{HI00} that appears in a construction
of a twistor space on $4\CP^2$ admitting a $\CC^*$-action,
whose algebraic dimension is two.

\medskip

As in Proposition \ref{prop:rs1},
the surface $S_0$ satisfies $K^2=0$.
Let $n$ be any integer greater than 4.
We then choose distinct $(n-4)$ points on the curve $C_0$,
any of which is different from the four points
$C_0\cap C_i,\,i\in\{1,2,3,4\}$, and  let 
\begin{align}\label{epsilon}
\epsilon:S\lra S_0
\end{align}
be the blowup at the $(n-4)$ points on $C_0$ and 
the $(n-4)$ points on $\ol C_0$ which are the real conjugations of the 
$(n-4)$ points.
Evidently the real structure and the $\CC^*$-action 
on $S_0$ naturally lift on the surface $S$.
Obviously $S$ is uniquely determined from $S_0$ by the choice
of $(n-4)$ points on $C_0$.
For the pluri-anticanonical systems of this surface,
we have

\begin{proposition}\label{prop:rs2}
For the above rational surface $S$, we have:
\begin{enumerate}
\item[\em (i)]
$K_S^2=8-2n$.
\item[\em (ii)]
For any $m\ge 1$, we have:
\begin{equation*}
h^0(mK_S\inv) =
\begin{cases}
0 & {\text{if $m$ is odd,}}\\
1 & {\text{if $m$ is even}}.
\end{cases}
\end{equation*}
In particular $\kappa\inv(S)=0$ for the
anti-Kodaira dimension of $S$.
\item[\em (iii)]
The single member of the system $|2K_S\inv|$ is 
concretely given by 
\begin{align}\label{isom0}
\Big(
2C_0 + \sum_{1\le i\le 4}C_i \Big)
+
\Big(
2\ol C_0 + \sum_{1\le i\le 4}\ol C_i \Big).
\end{align}
\end{enumerate}
\end{proposition}

\proof The item (i) is of course obvious.
For (ii), we recall from Proposition \ref{prop:rs1} the relation
\begin{align}\label{ac1}
K_0\inv \simeq f_0^*\ms O_{\CP^1}(1)
\end{align}
for the anti-canonical system on $S_0$.
Let 
$$
C_5, C_6,\dots, C_n
$$
be the exceptional curves over the blown-up points
on $C_0$, and $\ol C_5,\ol C_6,\cdots,\ol C_n$ the exceptional
curves over the real conjugate points over $\ol C_0$.
Let 
$$
f:S\stackrel{\epsilon}{\lra} S_0 
\stackrel{f_0}{\lra} \CP^1
$$
be the composition, which is also an elliptic fibration.
In the following, for simplicity of notation, we write
$$
C_{1,2,3,4}:=\sum_{1\le i\le 4} C_i\quad{\text{and}}\quad
C_{5,6,\dots,n}:=\sum_{5\le i\le n} C_i.
$$
Then from \eqref{ac1} we obtain
\begin{align}\label{ac2}
K_S\inv &\simeq \epsilon^* K_0\inv - 
(C_{5,6,\dots,n} + \ol C_{5,6,\dots,n})
\simeq \epsilon^*f_0^*\ms O(1) - (C_{5,6,\dots,n} + \ol C_{5,6,\dots,n})\notag\\
 &\simeq
f^*\ms O(1)-(C_{5,6,\dots,n} + \ol C_{5,6,\dots,n}).
\end{align}
Because $C_i$ and $\ol C_i$ belong to different fibers
of $f$, this in particular implies $H^0(K_S\inv) = 0$.
Further recalling that the components $C_0$ and $\ol C_0$ are included by multiplicity two
as in \eqref{I_0*}, we have
\begin{align}\label{ac4}
f^*\ms O(1) 
\simeq
2C_0 + C_{1,2,3,4} + 2C_{5,6,\dots,n},
\quad
f^*\ms O(1) 
\simeq
2\ol C_0 + \ol C_{1,2,3,4} + 2\ol C_{5,6,\dots,n}.
\end{align}
Therefore from \eqref{ac2} we have
\begin{align}\label{ac6}
2K_S\inv &\simeq f^*\ms O(2)- 2
(C_{5,6,\dots,n}+\ol C_{5,6,\dots,n})\notag\\
&\simeq
(
f^*\ms O(1) - 2C_{5,6,\dots,n}
)
+
(
f^*\ms O(1) - 2\ol C_{5,6,\dots,n}
)
\notag\\
&\simeq
(2C_0 + C_{1,2,3,4})
+
(
2\ol C_0 + \ol C_{1,2,3,4}).
\end{align}
Thus the last curve belongs to the system $|2K_S\inv|$,
and therefore
$h^0(2K_S\inv)\ge 1$,
implying $h^0(2mK_S\inv)\ge 1$ for any $m>0$.

Next we show $\kappa\inv(S) = 0$.
For this it is enough to see that 
the divisor \eqref{ac6} itself is precisely the `negative part'
of the Zariski decomposition \cite{Z62} of the divisor.
For this it suffices to verify that the intersection matrix
of the divisor
$
2C_0 + C_{1,2,3,4}
$
is negative definite.
From elementary calculation, we readily see that 
the eigenvalues of the intersection matrix
$(C_i, C_j)_{0\le i,j\le 4}$
are all negative.
Hence the intersection matrix is actually negative definite,
and thus we obtain $\kappa\inv (S) =0$.
This means that $h^0(2mK_S\inv) = 1$ for any $m>0$,
since otherwise we have $\kappa\inv(S)\ge 1$.

It remains to show that 
$H^0((2m-1)K_S\inv) = 0$ for any $m\ge 1$.
The case $m=1$ being already shown,  suppose $m>1$.
Since $K_S\inv.C_0=C_0^2 + 2 = (2-n) + 2=4-n$ which is negative as $n>4$,
the curve $C_0$ is a fixed component of $|(2m-1)K_S\inv|$.
As $K_S\inv.C_i = 0$ for any $i\in\{1,2,3,4\}$,
and $C_i$ and $C_0$ intersect,
this means that $C_i$ is also a fixed component
for any $i\in\{1,2,3,4\}$.
By reality, the same is true for the 
conjugate curves $\ol C_i$,
$0\le i\le 4$.
Hence we have an isomorphism
\begin{align}\label{ac37}
|(2m-1)K_S\inv| \simeq \big|(2m-1)K_S\inv - 
C_0 - C_{1,2,3,4} 
-\ol C_0 -\ol C_{1,2,3,4} \big|.
\end{align}
Moreover, the intersection number of the ingredient of 
the right-hand side with the curve $C_0$ is readily computed to be
$
2\{(1-m)n + (4m-5)\}.
$
Since $n\ge 5$ and $m>1$, we have $(1-m)n + (4m-5)\le 5(1-m)+(4m-5)
= -m<0$.
This implies that the right-hand side of \eqref{ac37} still has
$C_0$, and therefore $\ol C_0$ also, as  fixed components.
Thus the unique member \eqref{ac6} of $|2K\inv_S|$ is a fixed component
of the system $|(2m-1)K\inv_S|$ if $m>1$.
Hence we obtain an isomorphism
$$
|(2m-1)K_S\inv| \simeq |(2m-3)K_S\inv|
$$
for any $m>1$. As $|K_S\inv|=\emptyset$, this implies $|(2m-1)K_S\inv|=\emptyset$ 
for any $m\ge 1$, as desired.
\proofend


\section{Geometry of a twistor space which contains the rational surface}
\label{s:Z1}
Let $n>4$ and $g$ be a self-dual Riemannian metric 
on $n\CP^2$ whose scalar curvature is positive.
Let $Z$ be the associated twistor space.
In the sequel we write $F$ for the holomorphic line bundle $K^{-1/2}_Z$ which 
makes a natural sense for a twistor space associated to 
a self-dual metric \cite{AHS}.
In this section we investigate  structure of $Z$
which has the rational surface $S$ constructed in the last section as
a real member of the linear system $|F|$.
(The existence of such a twistor space will be shown in
Section 5.)

\subsection{Determination of the algebraic dimension of $Z$}
Let $n>4$ and $S$ be any one of the rational surfaces constructed
in the last section, satisfying $K^2 = 8-2n$.
Suppose the twistor space $Z$ has the surface $S$ as a real 
member of the system $|F|$.
We also suppose that the $\CC^*$-action on $S$ 
extends to the whole of $Z$ in a way that it is compatible with 
the real structure on $Z$.
This means that the corresponding self-dual structure on $n\CP^1$ 
have
a non-trivial $S^1$-action, and its natural lift
on the twistor space $Z$ has the given $\CC^*$-action
as the complexification
when restricted to the divisor $S$.
Then we have the following properties on
linear systems on $Z$:

\begin{proposition}\label{prop:Z1}
Suppose $n>4$, and let $Z$ and $S$ be a twistor space on $n\CP^2$
and the divisor respectively as above.
Then we have:
\begin{enumerate}
\item[\em (i)] The system $|F|$ consists of a single member $S$.
\item[\em (ii)] The system $|2F|$ ($=|K_Z\inv|$) is a pencil,
and the associated meromorphic map gives an algebraic reduction of $Z$.
In particular, $a(Z)=1$ holds for the algebraic dimension of $Z$.
\end{enumerate}
\end{proposition}

The meromorphic map $Z\to \CP^1$ in (ii) has a non-empty
indeterminacy locus. This will be investigated in the next
subsection.

The assertion (i) can be readily shown by a standard argument 
using Proposition \ref{prop:rs2}, and we omit a proof.
From (i), it follows $a(Z)\le 1$, since
we have $a(Z)\le 1 + \kappa\inv(S)$ in general
\cite{C91-2}, and $\kappa\inv(S) = 0$ by Proposition \ref{prop:rs2}.
The rest of this subsection is devoted to proving 
$a(Z)>0$.
We first show a transformation formula for the Euler characteristic of a line bundle 
over a compact complex threefold under the blowup 
with a 1-dimensional center in general:

\begin{proposition}\label{prop:RR}
Let $X$ be a compact complex threefold and
$L$ a holomorphic line bundle over $X$.
Let $C$ be a smooth curve on $X$,
$\mu: X_1\to X$  the blowup at $C$,
and $E\subset X_1$ the exceptional divisor.
Then if we write $L_1:=\mu^*L - E$, for the Euler characteristic, we have
\begin{align}\label{RR1}
\chi(L_1) = \chi(L) - L. C - 1+g_C,
\end{align}
where $g_C$ is the genus of the curve $C$.
\end{proposition}

\proof
By the Riemann-Roch formula for a threefold, we have
\begin{align}\label{RR2}
\chi(L_1) &= \frac 16 L_1^3 + \frac 14 L_1^2.\,c_1(X_1)
+ \frac1{12} L_1.\,\big(c_1(X_1)^2 + c_2(X_1)\big) + \frac 1{24}c_1(X_1).\,c_2(X_1).
\end{align}
To compute the right-hand side,
we first note that, for any line bundles $F_1$ and $F_2$ over $X$,
we have
\begin{align}\label{zero1}
\mu^*F_1.\mu^*F_2.E = \big(\mu^*F_1|_E,\mu^*F_2|_E\big)_E
= \big(\mu^*(F_1|_C),\mu^*(F_2|_C)\big)_E=0.
\end{align}
For $L_1^3$,  we readily have
\begin{align}\label{int1}
E^3 &= (E|_E)^2= - \deg N_{C/X},
\\
\label{int2}
(\mu^*L)^2.\,E &= 0\quad({\text{by }}\eqref{zero1}),\\
\label{int3}
(\mu^*L).\,E^2 &= (\mu^*L|_E, E|_E)_E = 
(\mu^*(L|_C), E|_E)_E
=
-L.C.
\end{align}
Therefore we obtain
\begin{align}\label{int3.5}
L_1^3 = L^3 - 3(\mu^*L)^2.\,E+ 3(\mu^*L).\,E^2
- E^3 = L^3 - 3L.C + \deg N_{C/X},
\end{align}
where we used \eqref{int1}--\eqref{int3} for the last equality.

Next for the second term in the right-hand side of \eqref{RR2}, we have
\begin{align}\label{int6}
L_1^2.\,c_1(X_1) &= \{(\mu^*L)^2 - 2\mu^*L.\,E +E^2\}.
(\mu^*c_1(X) - E)\notag\\
&=L^2.\,c_1(X) 
- 2L.C + E^2.\,\mu^*c_1(X) - E^3
\quad({\text{by }}\eqref{zero1}, \eqref{int3}).
\end{align}
Moreover we have
\begin{align}\label{int4}
E^2.\mu^*c_1(X) = \big(\mu^*K_X\inv|_E, E|_E\big)_E = 
\big(\mu^*(K_X\inv|_C), E|_E\big)_E = -K_X\inv.C
= 2g_C-2 -\deg N_{C/X},
\end{align}
where the last equality is from adjunction formula.
Hence, with the aid of \eqref{int1}, 
we obtain from \eqref{int6}
\begin{align}\label{int8}
L_1^2.\,c_1(X_1) &= L^2.\,c_1(X) - 2L.C +2g_C - 2.
\end{align}
Next, for the third term of \eqref{RR2}, we first have
\begin{align}\label{int9}
L_1.\,c_1(X_1)^2 &= (\mu^*L - E).(\mu^*c_1(X)-E)^2\notag\\
&= L.\,c_1(X)^2 + \mu^*L.\,E^2 + 2\mu^*c_1(X).\,E^2 - E^3
\quad({\text{by }} \eqref{zero1})\notag\\
&=L.\,c_1(X)^2 -L.C + 4(g_C-1)-\deg N_{C/X}.
\quad({\text{by }} \eqref{int3},\eqref{int4},\eqref{int1}).
\end{align}
On the other hand, for the transformation law for the second Chern class
of a threefold,
we have by \cite[pp.609--610, Lemma]{GH}
$$
c_2(X_1) = \mu^*(c_2(X) + [C]) - (\mu^* c_1(X)).E,
$$
where $[C]$ denotes the cohomology class in $H^4(X,\ZZ)$
represented by the 2-cycle $C$.
Therefore we obtain
\begin{align*}
L_1. c_2(X_1) &= (\mu^*L - E). \big\{\mu^*(c_2(X) + [C]) - (\mu^* c_1(X)).E\big\}\\
&= L.c_2(X) + L.C - E.\mu^*(c_2(X)+[C])
+ \mu^*c_1(X).E^2\quad({\text{by }}\eqref{zero1})\\
&= L.c_2(X) + L.C - E.\mu^*(c_2(X)+[C])
+2g_C-2-\deg N_{C/X} \quad({\text{by }} \eqref{int4}).
\end{align*}
For the third term, we have
$$
E.\mu^*(c_2(X)+[C]) = \big(E|_E,\,\mu^*(c_2(X) + [C])|_E\big)_E
= \big(E|_E,\, \mu^*((c_2(X) + [C])|_C\big)_E = 0,
$$
because $c_2(X) + [C]\in H^4(X,\ZZ)$ and hence
$(c_2(X) + [C])|_C=0$.
Therefore we obtain
\begin{align}\label{int10}
L_1. c_2(X_1) = L.c_2(X) + L.C 
+2g_C-2-\deg N_{C/X}. 
\end{align}
Substituting \eqref{int3.5},
\eqref{int8},
\eqref{int9} and \eqref{int10}  to the Riemann-Roch formula \eqref{RR2},
and noting that the product $c_1.\,c_2$ is a birational invariant,
we obtain the required formula \eqref{RR1}.
\proofend

\medskip

For proving the assertion (ii) in Proposition \ref{prop:Z1},
 we next recall that, on the divisor $S$,
 there are $(n-4)$ numbers of $(-1)$-curves
 $C_5,C_6,\dots,C_{n}$.
(See the construction in the last section.)
Since the restriction of the twistor projection 
$\pi:
Z\to n\CP^2$ to the divisor $S$ is necessarily diffeomorphic 
around the curves
$C_i$ and $\ol C_i$ for any $i\in \{5,6,\dots,n\}$
in the orientation-reversing way,
the image $\pi(C_i)=\pi(\ol C_i)$ is a
sphere in $n\CP^2$ whose self-intersection number is 1.
Moreover, these spheres are disjoint since
$\pi$ is degree two on $S$.
Let $\xi_i\in H^2(n\CP^2,\ZZ)$, $i\in\{5,6,\dots,n\}$, be 
the cohomology class represented by the sphere $\pi(C_i)$.
These form a part of orthonormal basis of $H^2(n\CP^2,\ZZ)\simeq \ZZ^n$.
We put $\aaa_i:=\pi^*\xi_i\in H^2(Z,\ZZ)$.
Since $H^1(\ms O_Z)=H^2(\ms O_Z)=0$, we have 
an isomorphism $\Pic\,Z\simeq H^2(Z,\ZZ)$, and 
we write $\ms O_Z(\aaa_i)$ when we regard the element
$\aaa_i\in H^2(Z,\ZZ)$ as an element in $\Pic\,Z$.
We then have
\begin{align}\label{rest20}
\ms O_Z(\aaa_i)|_S \simeq\ms O_S(C_i-\ol C_i).
\end{align}

Then we have the following key vanishing property:

\begin{lemma}\label{lemma:van1}
Let $Z$ and $S$ be as in Proposition \ref{prop:Z1},
and $C_0\subset S$ the rational curve on $S$ with 
self-intersection $(2-n)$ in $S$.
(See \eqref{I_0*} for the definition of $C_0$.)
Then we have
$$
H^1(\ms O_Z(-\aaa_5 - \aaa_6 - \dots - \aaa_n)
\otimes\ms I_{C_0})=0
$$
for the ideal sheaf $\ms I_{C_0}$ in $\ms O_Z$.
\end{lemma}

\proof
In this proof we write $\aaa_{5,6,\dots, n}$ 
and $C_{5,6,\dots,n}$ for
 the sum $\sum_{5\le i\le n}\aaa_i$ 
and $\sum_{5\le i\le n}C_i$ respectively, for simplicity.
From the inclusion $C_0\subset S\subset Z$ we have
an exact sequence
\begin{align*}
0 \lras \ms I_S \lras \ms I_{C_0}\lras \ms O_S(-C_0)\lras 0,
\end{align*}
where the two ideal sheaves are taken in $\ms O_Z$.
As $S\in |F|$ we have $\ms I_S\simeq -F$.
By taking a tensor product with the sheaf
$\ms O_Z(-\aaa_{5,6,\dots, n})$,
noting
\begin{align*}
 \ms O_Z(-\aaa_{5,6,\dots, n})\otimes
 \ms O_S(-C_0) \simeq
 \ms O_S(-C_0 - C_{5,6,\dots, n} + \ol C_{5,6,\dots, n})
\end{align*}
by \eqref{rest20} from the choice of the class $\aaa_i$, we obtain
an exact sequence
\begin{align}\label{ses5}
0 \to \ms O_Z(-\aaa_{5,6,\dots, n})
\otimes (-F) \to 
\ms O_Z(-\aaa_{5,6,\dots, n})\otimes\ms I_{C_0}\to
 \ms O_S(-C_0 - C_{5,6,\dots, n} + \ol C_{5,6,\dots, n})
 \to 0.
\end{align}
Now from the Hitchin's vanishing theorem \cite{Hi80},
noting that the line bundle
$\ms O_Z(\pm \aaa_{5,6,\dots, n})$ are real and of degree $0$,
 we have
\begin{align*}
H^1(\ms O_Z(-\aaa_{5,6,\dots, n})\otimes (-F))&=0,\\
H^2(\ms O_Z(-\aaa_{5,6,\dots, n})\otimes (-F))
&\simeq
H^1(\ms O_Z(\aaa_{5,6,\dots, n})\otimes (-F))^*=0.
\end{align*}
Therefore from the cohomology exact sequence of \eqref{ses5}, 
we obtain an isomorphism
\begin{align}\label{isom1}
H^1(\ms O_Z(-\aaa_{5,6,\dots, n})\otimes \ms I_{C_0})
\simeq
H^1(\ms O_S(-C_0 - C_{5,6,\dots, n} + \ol C_{5,6,\dots, n})).
\end{align}
In order to compute the right-hand side,
we consider an obvious exact sequence
\begin{align}\label{ses6}
0\lra 
\ms O_S(-C_0 - C_{5,6,\dots, n} + \ol C_{5,6,\dots, n})
\lra \ms O_S(\ol C_{5,6,\dots, n})
\lra \ms O_{C_0\cup C_5\cup C_6 \cup \dots\cup C_n}
\lra 0.
\end{align}
Here we have used the fact that 
the curves $\ol C_5,\ol C_6,\dots, \ol C_n$ are
disjoint from the union $C_0\cup C_5\cup C_6 \cup \dots\cup C_n$.
Then noting that 
$C_0\cup C_5\cup C_6 \cup \dots\cup C_n$ is simply connected
from their configuration,
we have $H^1(\ms O_{C_0\cup C_5\cup C_6 \cup \dots\cup C_n}
)=0$.
Moreover, we clearly have 
$H^0(\ms O_S(-C_0 - C_{5,6,\dots, n} + \ol C_{5,6,\dots, n}))=0$.
Therefore from the connectedness of 
$C_0\cup C_5\cup C_6 \cup \dots\cup C_n$,
 the cohomology exact sequence of \eqref{ses6} gives
an isomorphism 
\begin{align}\label{isom2}
H^1(\ms O_S(-C_0 - C_{5,6,\dots, n} + \ol C_{5,6,\dots, n}))
\simeq
H^1(\ms O_S(\ol C_{5,6,\dots, n})).
\end{align}
Moreover, since $\ol C_5, \ol C_6,
\dots, \ol C_n$ are $(-1)$-curves and the surface $S$ is rational,
we readily obtain $H^1(\ms O_S(\ol C_{5,6,\dots, n}))=0$.
Hence from \eqref{isom1} we obtain 
the required vanishing result.
\proofend

\medskip
For the purpose of proving the assertion (ii) of 
Proposition \ref{prop:Z1},
let
\begin{align}\label{bu1}
\mu_1:Z_1\to Z
\end{align}
be the blowup of the twistor space $Z$ at the curves 
$C_0\sqcup \ol C_0$, and 
let $E_0$ and $\ol E_0$ be the exceptional divisors over
$C_0$ and $\ol C_0$ respectively.
From the inclusion $C_0\subset S\subset Z$, recalling that 
$(C_0)^2_S = 2-n$, for the normal bundle
of $C_0$ in $Z$, we readily have
\begin{align}\label{nb1}
N_{C_0/Z}\simeq \ms O(3-n)\oplus \ms O(3-n)
\quad{\text{or}}\quad
\ms O(2-n)\oplus \ms O(4-n).
\end{align}
Accordingly, we have
\begin{align}\label{exc1}
E_0\simeq \CP^1\times \CP^1
\quad{\text{or}}\quad \FF_2.
\end{align}
Then by utilizing the transformation formula in Proposition \ref{prop:RR}, we can prove
\begin{lemma}\label{lemma:chi0}
Let $\mu_1:Z_1\to Z$ be as above.
Then we have
\begin{align}\label{chi0}
\chi\Big(Z_1,
\ms O_{Z_1}(-E_0)\otimes\mu_1^*\ms O_Z\big(
\aaa_5 + \aaa_6 + \dots + \aaa_n
\big)\Big)=0
\end{align}
for the Euler characteristic of the invertible sheaf on $Z_1$.
\end{lemma}

\proof
We decompose the blowup $\mu_1:Z_1\to Z$ as 
$$
Z_1\stackrel{\nu_2}{\lra} Z_{1/2}
\stackrel{\nu_1}{\lra} Z,
$$
where $\nu_1$ is the blowup at $C_0$ and $\nu_2$ is the blowup at $\ol C_0$.
(We use this strange notation only in this proof.)
Then we have, with the same notation as in the proof of
Lemma \ref{lemma:van1}, noting $\ol C_0\cap E_0=\emptyset$ in $Z_{1/2}$,
\begin{align}\label{isom3}
\mu_1^*\ms O_Z(-\aaa_{5,6,\dots,n})- E_0
\simeq
\nu_2^*
\Big(
\nu_1^*\ms O_Z(-\aaa_{5,6,\dots,n}) - E_0
\Big).
\end{align}
Moreover, since $\nu_2$ is birational, we have,
for any coherent sheaf $L$ on $Z_{1/2}$,
\begin{align*}
H^q(Z_1,\nu_2^*L)\simeq H^q(Z_{1/2},L)
\quad{\text{for any }} q\ge 0.
\end{align*}
Therefore from \eqref{isom3} we have
\begin{align*}
H^q
\Big(
Z_1,\,
\mu_1^*\ms O_Z
\big(-\aaa_{5,6,\dots,n}
\big)-E_0
\Big)
\simeq
H^q
\Big(
Z_{1/2},\,
\nu_1^*\ms O_Z
\big(-\aaa_{5,6,\dots,n}
\big)-E_0
\Big)
\quad{\text{for any }} q\ge 0.
\end{align*}
Hence we have
\begin{align}\label{isom4}
\chi
\Big(
Z_1,\,
\mu_1^*\ms O_Z
\big(-\aaa_{5,6,\dots,n}
\big)-E_0
\Big)
=
\chi
\Big(
Z_{1/2},\,
\nu_1^*\ms O_Z
\big(-\aaa_{5,6,\dots,n}
\big)-E_0
\Big).
\end{align}
Now we are going to compute the right-hand side 
by using Proposition \ref{prop:RR}
on putting $L=\ms O_Z
\big(-\aaa_{5,6,\dots,n}
\big)$, $C=C_0$, and $\mu=\nu_1$.
For this, we first compute as
\begin{align}
L.\,C_0 &= 
\big(
\ms O_Z
\big(
-\aaa_{5,6,\dots,n}
\big)|_S, \,C_0
\big)_S\notag\\
&= \big(
-C_{5,6,\dots,n} + 
\ol C_{5,6,\dots,n},\,C_0 
\big)_S \quad({\text{by \eqref{rest20}}})
\notag\\
&=-\big(
C_{5,6,\dots,n} ,\,C_0 
\big)_S \quad({\text{as }} \ol C_i\cap C_0 =\emptyset \,\,\,
{\text{for any }} i\in\{5,6,\dots,n\})\notag\\
&= 4-n \quad({\text{as }}  (C_i,\, C_0)_S =1 \,\,\,
{\text{for any }}i\in\{5,6,\dots,n\}).\label{int11}
\end{align}
Therefore by the formula in Proposition \ref{prop:RR}, we obtain,
recalling that $C_0$ is rational,
\begin{align}\label{chi1}
\chi
\Big(
Z_{1/2},\,
\nu_1^*\ms O_Z
\big(-\aaa_{5,6,\dots,n}
\big)-E_0
\Big)
&=
\chi
\Big(
Z,\,
\ms O_Z
\big(-\aaa_{5,6,\dots,n}
\big)
\Big) - (L,C_0)_Z -1.
\end{align}
The first term of right-hand side may be computed by Riemann-Roch formula
as follows:
\begin{align*}
\chi
\Big(
Z,\,
\ms O_Z
\big(-\aaa_{5,6,\dots,n}\big)
\Big)
&=
\frac16
\big(-\aaa_{5,6,\dots,n}\big)^3
+ \frac 14
\big(-\aaa_{5,6,\dots,n}\big)^2.\,c_1(Z)\\
&\qquad+\frac 1{12}
\big(-\aaa_{5,6,\dots,n}\big).\,
\big(c_1(Z)^2+c_2(Z)\big) + \frac1{24}c_1(Z)c_2(Z)\\
&=
\frac14 (-4)(n-4) + 1= 5-n,
\end{align*}
where in the second equality we used that 
\begin{align*}
\aaa_i.\,\aaa_j.\,\aaa_k&=0\quad{\text{for any }} i,j,k
\quad({\text{since $\aaa_i$-s are lifts from the 4-manifold}}),
\\
\aaa_i^2.\,c_1(Z) &= -4\quad{\text{for any }}  i\quad({\text{adjunction formula}}),\\
\aaa_i.\,c_1^2(Z) &= \aaa_i.\,c_2(Z) = 0
\quad({\text{as $c_1^2(Z)$ and $c_2(Z)$ are lifts 
from the 4-manifold \cite{Hi81}}}),\\
c_1(Z).\,c_2(Z)& = 12(\chi -\tau)= 24
\quad({\text{Hitchin \cite[p.135]{Hi81}}}).
\end{align*}
Thus from \eqref{chi1} and \eqref{int11}, we finally obtain
\begin{align*}
\chi
\Big(
Z_{1/2},\,
\nu_1^*\ms O_Z
\big(-\aaa_{5,6,\dots,n}
\big)-E_0
\Big) = (5-n) - (4-n) - 1 = 0.
\end{align*}
From \eqref{isom4}, this means the assertion \eqref{chi0}.
\proofend

\medskip
We are still preparing for the proof of Proposition \ref{prop:Z1} (ii).
We next show the following:

\begin{lemma}\label{lemma:isom2}
Under the  situation in Lemma \ref{lemma:chi0}, 
let $S_1\subset Z_1$ be the strict transform of $S$ into $Z_1$.
Then we have an isomorphism
\begin{align}\label{isom5}
\big(\mu_1^*
\big(
F-\aaa_{5,6,\dots,n}
\big)-2E_0\big)|_{S_1}
\simeq
\ms O_{S_1}(
C_1 + C_2 + C_3 + C_4
),
\end{align}
where we are identifying the curve $C_i$ in $S$ with 
the strict transform into $S_1$.
Moreover, we have an isomorphism
\begin{align}\label{isom6}
\big(\mu_1^*
\big(
F-\aaa_{5,6,\dots,n}
\big)-2E_0\big)|_{\ol E_0}
\simeq\ms O_{\ol E_0}.
\end{align}

\end{lemma}

\proof
Since the center $C_0\sqcup \ol C_0$ for the blowup $\mu_1$
is contained in the  divisor $S$,
the intersections $S_1\cap \ol E_0$
and $S_1\cap \ol E_0$ 
are naturally isomorphic to $ C_0$ and $\ol C_0$
under $\mu_1$ respectively.
Hence for the first restriction, we calculate,
under the identification $S_1\simeq S$ given by $\mu_1$,
\begin{align}
\big\{
\mu_1^*
\big(
F-\aaa_{5,6,\dots,n}
\big)-2E_0\big\}|_{S_1}
&\simeq
(F-\aaa_{5,6,\dots,n})|_S - 2C_0\notag\\
&\simeq \big(K_S\inv - C_{5,6,\dots,n} 
+\ol C_{5,6,\dots,n} \big) - 2C_0\notag\\
&\simeq
\big\{(
f^*\ms O(1) - C_{5,6,\dots,n} 
-\ol C_{5,6,\dots,n})
- C_{5,6,\dots,n} 
+\ol C_{5,6,\dots,n} 
\big\}-2C_0\notag\\
&\simeq f^*\ms O(1) - 2 C_{5,6,\dots,n} -2C_0\notag\\
&\simeq 
\big(
2C_0 + C_{1,2,3,4} + 2C_{5,6,\dots,n} 
\big)- 2 C_{5,6,\dots,n} -2C_0
\quad({\text{by }}  \eqref{ac4})\notag\\
&\simeq C_{1,2,3,4},\label{isom6.5}
\end{align}
where in the third isomorphism we have used the isomorphism
$K_S\inv\simeq f^*\ms O(1) - C_{5,6,\dots,n} 
-\ol C_{5,6,\dots,n}$, which is obtained in \eqref{ac2}.
This means the first isomorphism \eqref{isom5} in the lemma.
On the other hand, for the latter restriction, we have
\begin{align*}
\big\{
\mu_1^*
\big(
F-\aaa_{5,6,\dots,n}
\big)-2E_0\big\}|_{\ol E_0}
&\simeq
\big\{
\mu_1^*
\big(
F-\aaa_{5,6,\dots,n}
\big)
\big\}|_{\ol E_0}
\quad({\text{as }} E_0\cap \ol E_0=\emptyset)\\
&\simeq
\mu_1^*\big((F-\aaa_{5,6,\dots,n})|_{\ol C_0}\big).
\end{align*}
But since 
\begin{gather*}
(F,\ol C_0)_Z = (K_S\inv,\ol C_0)_S =(2-n)+2 = 4-n,\\
(\aaa_i,\ol C_0)_Z = (C_i - \ol C_i,\ol C_0)_S = 
-(\ol C_i,\ol C_0)_S = -1\quad{\text{for all}}\,\, i\in\{5,6,\dots,n\},
\end{gather*}
we obtain 
\begin{align*}
(F-\aaa_{5,6,\dots,n})|_{\ol C_0}
\simeq 
\ms O_{\ol C_0}((4-n) +(n-4))\simeq\ms O_{\ol C_0}.
\end{align*}
Hence we get
the second isomorphism \eqref{isom6}
in the lemma.
\proofend

\medskip
Now we are able to prove a key proposition,
which readily implies $a(Z)>0$:
\begin{proposition}\label{prop:nonvan1}
Let $Z$ and $S$ be as in Proposition \ref{prop:Z1}.
Then we have a non-vanishing
\begin{align}\label{isom7}
H^0
\Big(
Z, F\otimes\ms O_Z
(- \aaa_5 - \aaa_6 - \dots - \aaa_n)
\Big)\simeq\CC.
\end{align}
Therefore, the system
$\big|F - \aaa_5 - \dots - \aaa_n\big|$
consists of a single member.
\end{proposition}

\proof
We consider the restriction homomorphism of the line bundle
$\mu_1^*
\big(
F-\aaa_{5,6,\dots,n}
\big)-2E_0$ to the divisor $S_1\cup \ol E_0$.
The kernel of this homomorphism is computed
as, noting $\mu_1^*F \simeq \ms O_{Z_1}(S_1 + E_0 + \ol E_0)$,
\begin{align*}
\mu_1^*
\big(
F-\aaa_{5,6,\dots,n}
\big)-2E_0
- (S_1 + \ol E_0)
&\simeq
 S_1 + E_0 + \ol E_0
+ \mu_1^*\ms O_Z(-\aaa_{5,6,\dots,n})-2E_0- S_1 - \ol E_0\\
&\simeq
\mu_1^*\ms O_Z(-\aaa_{5,6,\dots,n})
-E_0.
\end{align*}
Therefore we get an exact sequence
\begin{align}\label{ses7}
0&\lras
 \mu_1^*\ms O_Z(-\aaa_{5,6,\dots,n})-E_0 \lras
\mu_1^*
\big(
F-\aaa_{5,6,\dots,n}
\big)-2E_0\notag\\
&\hspace{6cm}\lras
\big\{
\mu_1^*\big(
F-\aaa_{5,6,\dots,n}
\big)-2E_0\big\}|_{S_1\cup \ol E_0}
\lras 0.
\end{align}
For the section of the last non-trivial term,
we have, by \eqref{isom5} and \eqref{isom6},
recalling that $S_1\cap \ol E_0=\ol C_0$ transversally
and $(C_1\cup C_2\cup C_3\cup C_4)\cap \ol C_0=\emptyset$, 
\begin{align*}
H^0\Big(S_1\cup \ol E_0,\,
\mu_1^*\big(
F-\aaa_{5,6,\dots,n}
\big)-2E_0\Big) \simeq\CC.
\end{align*}
On the other hand, by the Leray spectral sequence
and Lemma \ref{lemma:van1},
we obtain
\begin{align*}
H^1
\Big(Z_1,\,
\mu_1^*\ms O_Z(-\aaa_{5,6,\dots,n})-E_0
\Big)
&\simeq
H^1\Big(Z,\,\ms O_Z(-\aaa_{5,6,\dots,n})\otimes\ms I_{C_0}\Big), 
\end{align*}
and this vanishes by Lemma \ref{lemma:van1}.
Moreover, $H^0(\mu_1^*\ms O_Z(-\aaa_{5,6,\dots,n})-E_0)=0$ clearly.
Hence from the cohomology exact sequence of \eqref{ses7}, 
we obtain the desired non-vanishing result.
\proofend

\medskip
In this and the next subsections
we denote by $X$ for the unique member of 
the system $|F - \aaa_5-\dots-\aaa_n|$.
Then $\ol X$ is  the unique member of 
the system $|F + \aaa_5+\dots+\aaa_n|$.

%

\medskip\noindent
{\em Proof for (ii) of Proposition \ref{prop:Z1}.}
We have $h^0(2F)\ge 2$, because
 $2S$ and $X+\ol X$ are linearly independent divisors
 belonging to $|2F|$.
For the reverse inequality,
the cohomology exact sequence of 
the exact sequence 
$0 \to F \to 2F \to 2K_S\inv \to 0$
 implies $h^0(2F)\le h^0(F) + h^0(2K_S\inv)$.
But the right-hand side is $1 + 1 = 2$ by
Propositions \ref{prop:rs2} and \ref{prop:Z1}.
Thus we obtain $h^0(2F) = 2$.
Hence $a(Z)\ge 1$.
Since we have $a(Z)\le 1$ as
explained right after Proposition \ref{prop:Z1},
 we obtain $a(Z)=1$.

Let $\Phi:Z\to\CP^1$ be the rational map 
induced by the pencil $|2F|$,
and let $Z'\to Z$ be an elimination of the indeterminacy 
locus of $\Phi$ which preserves the real structure.
(We will an explicit elimination later in the next subsection,
 but here we do not need it.)
If $\Phi$ is not an algebraic reduction,
we get a factorization
$Z'\to \Gamma\to \CP^1$ of 
the morphism $Z'\to Z\to \CP^1$ as the Stein factorization, where $\Gamma\to\CP^1$ is a $d:1$
($d>1$) map
from some curve $\Gamma$. 
Since $2F.l =4$ for a twistor line $l$,
only $d\in\{2,4\}$ can occur.
If $d=4$, a general fiber $D'$ of $Z'\to \Gamma$ satisfies
$D'.l'=1$ for the inverse image of a generic twistor line $l$.
This contradicts the absence of a divisor $D$ on $Z$ which
satisfies $D.l=1$ (as $|F|=\{S\}$).
Hence $d\neq 4$.
If $d=2$, then taking the real structure into account,
a general fiber of $Z'\to \CCC$ is a transform of a member of $|F|$ on $Z$.
This again contradict $|F|=\{S\}$.
Hence $\Phi$ itself is a (meromorphic) algebraic reduction
of $Z$.
\proofend

\medskip
We emphasize that the pencil $|2F|$ is generated by 
the non-reduced member $2S$ and a reducible member $X+\ol X$ as in the above proof.

\subsection
{Structure of members of the anti-canonical system of $Z$}
Let $Z\to n\CP^2$ and $S\in |F|$ be as in the previous subsection.
In this subsection we investigate  structure of 
the pencil $|2F|$.
In particular we determine structure of 
a general fiber of the algebraic reduction of $Z$.
We first investigate the unique member
$X$ of the system 
$|F - \sum_{5\le i\le n}\aaa_i|$ found in the last subsection.

\begin{proposition}\label{prop:X1}
Let $X$ be the divisor as above.
\begin{enumerate}
\item[\em (i)]
As a divisor on the surface $S\in |F|$, we have
\begin{align}\label{rest1}
X|_S = 2C_0 + \sum_{1\le i\le 4} C_i.
\end{align}
\item[\em (ii)] $X$ is irreducible, non-normal, and 
of multiplicity two along the curve $C_0$.
(At this stage we do not yet claim that 
the singularity is  ordinary double.)
\end{enumerate}

\end{proposition}

\proof
For (i), 
by a similar computation for obtaining \eqref{isom6.5},
we have
\begin{align*}
X|_S &\simeq 2C_0 + C_{1,2,3,4}.
\end{align*}
Further it is easy to see from intersection numbers that 
the system $|2C_0 + C_{1,2,3,4}|$ has no member other  than
$2C_0 + C_{1,2,3,4}$ itself.
Thus we obtain (i).

For the assertion (ii),
since the pencil $|2F|$ is generated by $2S$ and $X+\ol X$
as above,
by \eqref{rest1},
we have $C_0\cup\ol C_0\subset \Bs\,|2F|$.
Then while it is tempting to conclude that 
the divisor $X$ has double points along the curve $C_0$,
we need to exclude the possibility that 
$X$ is tangent to $S$ along $C_0$.
For this purpose,
as in the previous subsection,
let $\mu_1:Z_1\to Z$ be blowup at $C_0\sqcup\ol C_0$,
and $E_0\sqcup\ol E_0$ the exceptional divisor.
Then we  put 
\begin{align*}
L' := \mu_1^*(2F) - (E_0+\ol E_0).
\end{align*}
We have an identification
$|2F|\simeq |L'|.$
We now show that the system $|L'|$
has the divisor $E_0+\ol E_0$ as
 fixed components.

As in \eqref{nb1}, 
we have 
$N_{C_0/Z}\simeq \ms O(3-n)\oplus \ms O(3-n)$
or
$N_{C_0/Z}\simeq 
\ms O(2-n)\oplus \ms O(4-n)$,
and accordingly $E_0\simeq\FF_0$ or $E_0\simeq\FF_2$.
When the former is the case, writing $\ms O(0,1)$
for the fiber class of the projection
$E_0\to C_0$, we have
\begin{gather}\label{rest45}
E_0|_{E_0}\simeq \ms O(-1,3-n),\quad
(\mu_1^* F)|_{E_0}\simeq\mu_1^*(F|_{C_0})\simeq
\ms O(0,4-n).
\end{gather}
Hence we obtain 
\begin{align}\label{isom8}
L'|_{E_0}\simeq \ms O(0,8-2n)
\otimes\ms O(1,n-3)\simeq\ms O(1,5-n).
\end{align}
This means that, if $n>5$, 
the exceptional divisor $E_0$, and therefore $\ol E_0$ also,
are  fixed components of $|L'|$.
(The case $n=5$ will be considered later.)
When $N_{C_0/Z}\simeq 
\ms O(2-n)\oplus \ms O(4-n)$,
if we denote $A$ for the $(-2)$-section of 
the projection $E_0\simeq\FF_2\to C_0$,
and $\mathfrak f$ denotes its fiber class,
we have, in a similar way for obtaining \eqref{isom8},
\begin{align}\label{isom9}
L'|_{E_0}\simeq A -(n-6)\mathfrak f.
\end{align}
This again means $E_0\cup\ol E_0\subset\Bs\,|L'|$,
 if $n>6$.

If $E_0\simeq\FF_0$ and $n=5$, 
we can show $E_0\subset\Bs\,|L'|$ in the following way.
If this did not hold, by \eqref{isom8},
there would exist a divisor $Y'\in |L'|$
which satisfies $Y'|_{E_0}\in |\ms O(1,0)|$.
On the other hand,
if we still denote $C_i$ and $\ol C_i$
($i\in\{1,2,3,4\})$ for the strict transforms of 
the curves in $Z$ into $Z_1$,
the curve $C_i$ (and also $\ol C_i$)
is clearly a base curve of $|L'|$ by \eqref{rest1}.
In particular, $Y'$ has to include the four curves
$C_1,\dots,C_4$, and therefore
the point $C_i\cap E_0$ has to be on a $(1,0)$-curve for any $i\in\{1,2,3,4\}$.
On the other hand, if $S_1$ again denotes the strict transform
of the divisor $S$ into $Z_1$, we have by \eqref{rest45}
\begin{align}\label{isom10}
S_1|_{E_0}
\simeq \mu_1^* F - E_0 - \ol E_0|_{E_0}
\simeq \ms O(0,4-n)\otimes\ms O(1,n-3)
\simeq\ms O(1,1).
\end{align}
Moreover, since $\mu_1$ blows up the curve $C_0$ (and $\ol C_0$)
 in the  divisor $S$,
the intersection $S_1\cap E_0$ is isomorphic to the center
$C_0$, and therefore irreducible.
Furthermore, $S_1$ clearly contains the curves
$C_1,\dots,C_4$.
Therefore all the intersection points $E_0\cap C_i$
$(i\in\{1,2,3,4\})$ have to pass an irreducible 
$(1,1)$-curve on $E_0$.
This contradicts the above conclusion that the 4 points
$C_i\cap E_0$ are on a $(1,0)$-curve, which is derived
from the assumption  $E_0\not\subset \Bs\,|L'|$.
Therefore $E_0\subset \Bs\,|L'|$, as claimed.

Still we need to show that 
$E_0\subset \Bs\,|L'|$ holds also in the case
$E_0\simeq\FF_2$ and $n\in\{5,6\}$.
But this can be shown in a similar manner to the above argument,
so we omit the detail;
we just mention that, instead of \eqref{isom10}, we have
\begin{align}\label{isom10.5}
S_1|_{E_0}\simeq A + 2\mathfrak f.
\end{align}

Now since the pencil $|\mu_1^*(2F)|$ has the divisor
$2E_0+2\ol E_0$ as fixed components, any member
of the pencil $|2F|$ has multiplicity two along the curve
$C_0$. 
Hence so is the member $X+\ol X$.
But from \eqref{rest1} we have $\ol X\cap C_0=\emptyset$. 
Therefore
the divisor $X$ itself has  double points
along the curve $C_0$ at least.
But from \eqref{rest1}, $X$ can have at most double points
along $C_0$.
Therefore the multiplicity of $X$ along $C_0$ is two.
Finally, since $X\in |F-\aaa_{5,6,\dots,n}|$, if $X$ were reducible, we have $X = D + D'$,
for some degree-one divisors $D$ and $D'$.
Then $D+ \ol D$ is a member of $|F|$, and $D+\ol D\neq S$
as $S$ is irreducible.
This contradicts $h^0(F)=1$ in Proposition \ref{prop:Z1}.
Hence $X$ is irreducible.
\proofend

\medskip
Now the next proposition is obvious from \eqref{rest1}
if we recall that the pencil $|2F|$ is generated by 
the divisors $2S$ and $X+ \ol X$:

\begin{proposition}\label{prop:base35}
For the base locus of the anti-canonical system $|2F|$, we have
$$\Bs\,|2F| = \Big(\bigcup_{0\le i\le 4}C_i\Big)
\cup 
\Big(\bigcup_{0\le i\le 4}\ol C_i\Big).
$$
\end{proposition}

Therefore if  $\mu_1:Z_1\to Z$ is the blowup at $C_0\cup\ol C_0$ and
$E_0\cup\ol E_0$ is the exceptional divisor as before,
and if we put
\begin{align}\label{L1}
L_1: = \mu_1^*(2F) - 2(E_0+\ol E_0)
\end{align}
this time, then by Proposition \ref{prop:X1} (ii),  we have an isomorphism
$|2F|\simeq |L_1|$.
In order to determine structure of a general member of 
the original pencil $|2F|$, we look at the behavior of the pencil
$|L_1|$ when restricted to $E_1$:

\begin{lemma}\label{lemma:rest1}
\begin{enumerate}
\item[\em (i)]
We have an isomorphism
\begin{align}\label{isom11}
L_1|_{E_0}\simeq K_{E_0}\inv.
\end{align}
\item[\em(ii)]
The restriction map $H^0(Z_1,L_1)\to H^0(E_0,L_1)$ is injective.
Therefore the restriction of the pencil $|L_1|$ to $E_0$ 
remains to be a pencil.
\item[\em(iii)]
The last pencil on $E_1$ is generated by the twice of 
$S_1|_{E_0}\in |K^{-1/2}_{E_0}|$ which is smooth, and 
an anti-canonical curve which intersects
the curve $S_1|_{E_0}$ transversally at four points.
\item[\em (iv)]
A general member of the last pencil is a non-singular 
elliptic curve, which is $2:1$ over $C_0$ under
the projection $\mu_1|_{E_0}:E_0\to C_0$.
\end{enumerate}
\end{lemma}

\proof
Let $S_1\subset Z_1$ be the strict transform of $S$ into $Z_1$ as before.
Then since $S_1\in|\mu_1^*F - E_0 - \ol E_0|$, 
we have $2S_1\in|L_1|$. Hence
$
L_1|_{E_0}\simeq 2S_1|_{E_0}.
$
Then from \eqref{isom10} (in the case $E_0\simeq\FF_0$)
and \eqref{isom10.5} (in the case $E_0\simeq\FF_2$),
we obtain the desired isomorphism \eqref{isom11}.

Next, to show the injectivity of the restriction map,
let $X_1$ be the strict transform of the divisor 
$X$ into $Z_1$.
Then since the pencil $|L_1|$ is generated by 
$2S_1$ and $X_1+\ol X_1$, it suffices to 
show the two curves $2S_1|_{E_0}$ and
$(X_1+\ol X_1)|_{E_0}$ are different.
By \eqref{isom6.5}, we have
\begin{align}\label{bs2}
X_1|_{S_1}=
C_{1,2,3,4}.
\end{align}
This clearly means  $X_1|_{E_0}\neq S_1|_{E_0}$.
Moreover $\ol X_1\cap E_0=\emptyset$ by \eqref{rest1}
and hence $(X_1+\ol X_1)|_{E_0} = X_1|_{E_0}$.
Therefore $2S_1|_{E_0} \neq (X_1+\ol X_1)|_{E_0}$,
and we obtain (ii).

\begin{figure}
\includegraphics{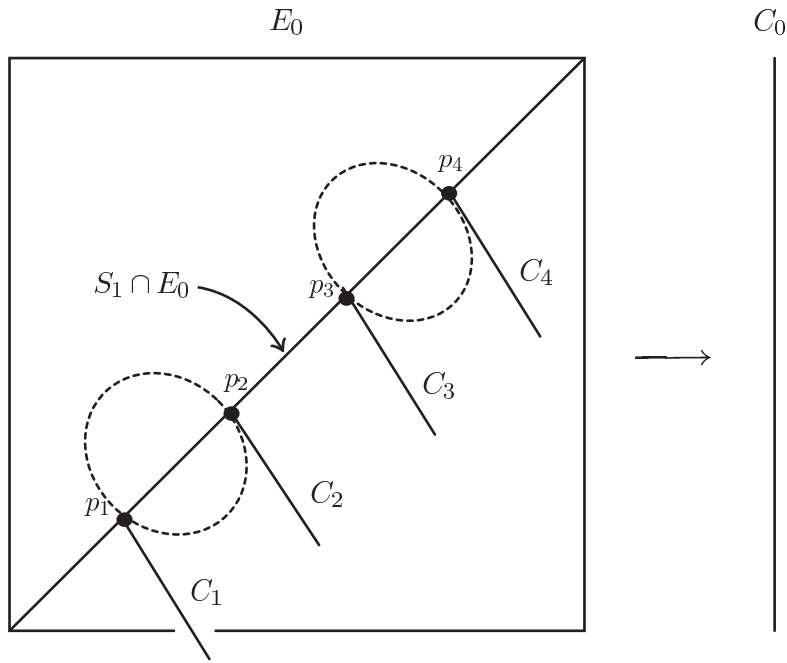}
\caption{Exceptional divisor of the blowup $\mu_1:Z_1\to Z_0$.
The curves $C_i$-s are transversal to $E_0$.
Dashed circles are a general member of the anti-canonical system on $E_0$.}
\label{fig:exe}
\end{figure}

For (iii), since the pencil $|L_1|$ is generated
by $2S_1$ and $X_1+\ol X_1$, 
 again writing $C_0$ for the intersection
$S_1\cap E_0$, it suffices to 
show that the restriction
$X_1|_{E_0}$, which is an anti-canonical by 
\eqref{isom11}, is a curve which intersects $C_0$
transversally at the 
four points $p_i:=C_i\cap E_0$, $i\in\{1,2,3,4\}$.
(See Figure \ref{fig:exe}.)
From \eqref{bs2}
we have $\{p_1,p_2,p_3,p_4\}\subset X_1$
and  $X_1\cap E_0\not\supset C_0$.
On the other hand we have 
$$
(X|_{E_0},C_0)_{E_0} = K_{E_0}\inv.\,K_{E_0}^{-1/2}
= (1/2)c_1^2(E_0) = (1/2)\cdot 8 = 4.
$$
This means $(X\cap E_0)\cap C_0 = \{p_1,p_2,p_3,p_4\}$
and each of the intersections is transverse.
Thus we obtain (iii).
This means that the base locus of 
the pencil is 
$\{p_1,p_2,p_3,p_4\}$.
Hence by Bertini's theorem, a general member of the pencil
can have singularities only at the 4 points.
But from the last transversality, this cannot happen.
Thus a general member of the pencil is a non-singular,
and it has to be an elliptic curve which is $2:1$ over
$C_0$, since it is an anti-canonical
curve. Thus we obtain (iv).
\proofend

\medskip
Later in the proof of Theorem \ref{thm:1} 
we will see that the intersection $X_1|_{E_0}$ is a smooth elliptic curve.

In order to reveal interesting geometric properties of the twistor space,
we determine the set of $\CC^*$-fixed points
on the twistor space.
For this, we recall that each of the curves $C_i$, $1\le i\le n$,
is $\CC^*$-invariant and there are exactly two fixed points
on each of these curves,
one of which is the intersection with the curve $C_0$.
Let $l_i$ be the twistor line that passes the other fixed point
on $C_i$.
These $n$ twistor lines are $\CC^*$-invariant.
Then we have:

\begin{lemma}\label{lemma:fl}
For each $i\in \{1,2,\dots, n\}$, let $l_i$ be the twistor line
as above.
Then the set of $\CC^*$-fixed points on the twistor space $Z$ consists of the  two curves $C_0$ and $\ol C_0$,
the eight points $l_i\cap(C_i\cup \ol C_i)$, $i\in \{1,2,3,4\}$,
and the $(n-4)$ twistor lines $l_i$, $i\in \{5,6,\dots, n\}$.
\end{lemma}

\proof
The divisor $S\in|F|$ is $\CC^*$-invariant,
and from the construction of $S$, 
the set of $\CC^*$-fixed points on $S$ consists
of the curve $C_0\cup\ol C_0$ and
one of the two fixed points on the curves 
$C_i$ and $\ol C_i$ for $1\le i\le n$.
From computations using local coordinates,
it is easy to see that
if $i\in\{1,2,3,4\}$, the twistor line $l_i$  has exactly two fixed points $l_i\cap(C_i\cup \ol C_i)$,
while
if $i\in \{5,6,\dots, n\}$, the twistor line $l_i$
is point-wise fixed by the $\CC^*$-action.
Then by Lefschetz fixed-point formula, noting $\chi(Z) = 2n+4$, it follows readily from this that these are all the fixed points of the $\CC^*$-action on $Z$.
\proofend

\medskip
We call these $(n-4)$ twistor lines $l_i$,
$i\in\{5,6,\dots,n\}$, as {\em $\CC^*$-fixed twistor lines}.

Now we can prove the main result of this subsection,
which in particular implies that the twistor spaces we have discussed
satisfy the property that a general fiber of its algebraic
reduction is birational to an elliptic ruled surface:

\begin{theorem}\label{thm:1}
Let $n>4$ and $Z$ be a twistor space on $n\CP^2$ with $\CC^*$-action,
which has the rational surface $S$ constructed in Section 2
as a real (single) member of $|F|$ as in Proposition \ref{prop:Z1}.
Then we have:
\begin{enumerate}
\item[\em (i)] 
A general member of the pencil $|2F|$ is birational to 
a ruled surface over an elliptic curve.
Moreover, the two rational curves $C_0$ and $\ol C_0$ 
are ordinary double curves of these surfaces.
\item[\em (ii)] 
The unique member $X$ of the system
$|F-\aaa_5-\dots-\aaa_n|$ in Proposition \ref{prop:X1}
is birational to a Hopf surface.
Moreover, the curve $C_0$ is an ordinary double curve of $X$.
\end{enumerate}
\end{theorem}

\proof 
The assertions about double curves in (i) and (ii) are now obvious from 
Lemma \ref{lemma:rest1} (iv).
For proving the remaining assertion in (i),
we first see that each member of the pencil
$|2F|$ is $\CC^*$-invariant.
For this we notice that each point of 
the image of $C_0$ under the twistor fibration
$Z\to n\CP^2$ is fixed by the $S^1$-action on $n\CP^2$
corresponding to the $\CC^*$-action on $Z$.
This means that under holomorphic coordinates
$(x,y,z)$ in a neighborhood of a point of $C_0$
for which $C_0=\{x=y=0\}$, the $S^1$-action 
is concretely given by $(x,y,z)\mapsto (tx,ty,z)$,
$t\in S^1=U(1)$.
Therefore the $\CC^*$-action induced on 
the exceptional divisor $E_0$ is trivial.
Since the restriction map
$H^0(Z_1,L_1)\to H^0(E_0,L_1)$ is injective as in 
Lemma \ref{lemma:rest1} and is  $\CC^*$-equivariant
from $\CC^*$-invariance of $C_0$,
this means that 
each member of the pencil
$|2F|$ is $\CC^*$-invariant.

Take a general member $Y$ of the pencil $|2F|$,
and let $Y_1\in |L_1|$ be the strict transform of $Y$ into 
the blowup $Z_1$.
Then $Y_1$ is also $\CC^*$-invariant, and
the two intersections $Y_1\cap E_0$ and $Y_1\cap\ol E_0$  are smooth elliptic curves on $E_0$ by Lemma \ref{lemma:rest1} (iv), whose points are $\CC^*$-fixed.
Therefore,
if $\tilde Y_1$ denotes an equivariant resolution of $Y_1$
(including the case $\tilde Y_1=Y_1$ if $Y_1$ is non-singular),
the strict transforms of the intersections $Y_1\cap E_0$ 
and $Y_1\cap\ol E_0$
into $\tilde Y_1$
are also $\CC^*$-fixed elliptic curves.
In particular, a non-singular minimal model of 
the surface $\tilde Y_1$ has a $\CC^*$-fixed point.
Therefore by the classification of compact complex surfaces
admitting a non-trivial $\CC^*$-action which has a fixed point obtained
in \cite{OW77} and \cite{Ha95},
$Y_1$ is birational to one of the following surfaces:
a rational surface, a ruled surface over a curve
of genus $\ge 1$, or a VII-surface.
But a non-trivial $\CC^*$-action on a rational surface cannot have an elliptic curve which is point-wise fixed.
Also,  a ruled surface over a curve
of genus $> 1$ does not have an elliptic curve.
Therefore $Y_1$ is birational to an elliptic ruled surface or
a VII-surface.
But by the classification of a non-trivial $\CC^*$-action on a VII-surface
with a fixed point obtained in \cite{Ha95},
such a $\CC^*$-action on a VII-surface cannot have two elliptic curves
which are point-wise fixed.
Hence $Y_1$ is birational to an elliptic ruled surface,
and we obtain (i) of the theorem.

For the assertion (ii), let $X_1$ be the strict transform
of $X$ into $Z_1$.
Then the intersection $X_1\cap E_0$ is a curve whose points
are fixed by the induced $\CC^*$-action.
Let $\nu:\tilde X_1\to X_1$ be a 
desingularization which preserves the $\CC^*$-action.
$(\nu_1$ might be an identity map.)
We first show that $\tilde X_1$ does not contain 
an irreducible curve $D$ which is point-wise fixed by the $\CC^*$-action
and which is different from any component
of the curve $\nu\inv(X_1\cap E_0$).
If the image $\mu_1(\nu(D))$ is a curve, then it must be one of the fixed twistor line
$l_i$-s from Lemma \ref{lemma:fl}.
But if this were actually the case,
we would have $l_i\subset X$ for some fixed twistor line $l_i$
($5\le i\le n$).
This contradicts
\begin{align}\label{rest11}
X|_S = 2C_0 + C_{1,2,3,4}
\end{align}
obtained in Proposition \ref{prop:X1},
 since $l_i$ is disjoint from the right-hand side
 of \eqref{rest11} if $i\in\{5,6,\dots,n\}$.
Thus we have $\mu_1(\nu(D))\neq l_i$ for any $i\in\{5,6,\dots,n\}$.
Hence $\mu_1(\nu(D))$ has to be a point (if $D$ would exist).
But if $\mu_1(\nu(D))$ is a point, then  the closure of a general orbit of the $\CC^*$-action on $X$
passes the point $\mu_1(\nu(D))$.
This cannot happen as is readily seen from the fact that 
the point on $n\CP^2$ over which the point $\mu_1(\nu(D))$ belongs
must be an isolated fixed point of the $S^1$-action corresponding to 
the $\CC^*$-action on $Z$.
Thus we have shown that the desingularization
$\tilde X_1$ does not contain 
an irreducible curve $D$ which is point-wise fixed by the $\CC^*$-action
and which is different from any component
of $\nu\inv(X_1\cap E_0$).

Next by using this we show that $\tilde X_1$ is a VII-surface.
From Lemma \ref{lemma:rest1} (iii),  the intersection 
$X_1\cap E_0$ is an anti-canonical curve on $E_0$.
Moreover, regardless of whether $E_0\simeq\FF_0$ or $\FF_2$,
it can be readily seen that the curve $X_1|_{E_0}$ does not have a non-reduced component.
This curve $X_1|_{E_0}$ is point-wise fixed by the $\CC^*$-action.
Hence if $X_1$ is birational to a ruled surface,
from the absence of the curve $D$ above,
there has to be a point on $X_1$ such that 
the closure of a general orbit on $X_1$ of the $\CC^*$-action
intersects.
This cannot  happen as before, and therefore 
$X_1$ cannot be a ruled surface.
Hence by the above cited result of \cite{Ha95},
$\tilde X_1$ is birational to a VII-surface.
Then a finer classification result  \cite[Klassifikationssatz]{Ha95} means that 
if a VII-surface admits a non-trivial $\CC^*$-action which has a fixed point,
then the surface is a Hopf surface or a parabolic Inoue surface
\cite{I74, E81}.
Both of these surfaces with $\CC^*$-action has
exactly two connected curves, one of which is a non-singular
elliptic curve which is point-wise fixed by the 
$\CC^*$-action.
The two possibilities are distinguished by another curve;
if the surface is Hopf, then the curve is a non-singular elliptic curve which is an orbit of the $\CC^*$-action,
and if the surface is parabolic Inoue, then 
the curve is a cycle of rational curves.

We show that the surface cannot be a parabolic Inoue surface.
For this, we first note that the intersection $X_1|_{E_0}$
has to be a smooth elliptic curve from the above description
of the $\CC^*$-action on VII-surfaces.
This in particular means that $X_1$ is smooth at
points on $ X_1|_{E_0}$.
So the self-intersection number in $X_0$ makes sense.
Recalling that $X_1\in|\mu_1^*(F-\aaa_{5,6,\dots,n})-2E_0|$,
it can be computed as
\begin{align}\label{100}
(E_0|_{X_1})^2 &= E_0.E_0.X_1=E_0.E_0.\{\mu_1^*(F-\aaa_{5,6,\dots,n})-2E_0\}.
\end{align}
Moreover we have
\begin{align}\label{101}
E_0.E_0.\mu_1^*F &= \big(\mu_1^*(F|_{C_0}), E_0|_{E_0}\big)_{E_0}\notag\\
&= -F.C_0 = -K_S\inv.C_0 = n-4,\\
E_0.E_0.\mu_1^*(\aaa_{5,6,\dots,n}) &= 
\big(\mu_1^*(\aaa_{5,6,\dots,n}|_{C_0}), E_0|_{E_0}\big)_{E_0}\notag\\
&= -(C_{5,6,\dots,n}-\ol C_{5,6,\dots,n},C_0)_S = -(n-4).
\end{align}
Further, in the case $E_0\simeq \FF_0$, we have
\begin{align}\label{103}
E_0^3 &= 
(E_0|_{E_0})^2\notag\\
&= (-1,3-n).(-1,3-n) = 2(n-3),\end{align}
and the same conclusion for the case $E_0\simeq\FF_2$.
Substituting \eqref{101}--\eqref{103} into \eqref{100},
we obtain 
\begin{align}\label{104}
(E_0|_{X_1})^2 &= (n-4) + (n-4) -4(n-3) \notag\\
&= 2(2-n).
\end{align}
Therefore, if $X_{\rm{min}}$ denotes the minimal model
of $X_0$, since the self-intersection number of the
elliptic curve in $X_{\rm{min}}$  corresponding to $E_0|_{X_1}$ is zero,
we can obtain the minimal model $X_{\rm{min}}$ from the desingularization
$\tilde X_1$ by blowing-up $(-1)$-curves that intersect
the elliptic curve $X_1|_{E_0}$ $2(n-2)$ times repeatedly.
In particular, $\tilde X_1$ has at least $2(n-2)$
isolated fixed points of the $\CC^*$-action.

We next show that the curves $C_1,C_2,C_3$ and $C_4$
(the strict transforms into $Z_1$)
are $(-1)$-curves on $X_1$.
For this, from the inclusions $C_i\subset S_1\subset X_1$,
we obtain an exact sequence
$0\to N_{C_i/S_1}\to N_{C_i/Z_1}\to N_{S_1/Z_1}|_{C_i}\to 0$.
Obviously we have $N_{C_i/S_1}\simeq N_{C_i/S}\simeq\ms O(-2)$.
Further, we have 
$N_{S_1/Z_1}\simeq N_{S/Z} - C_0 - \ol C_0
\simeq K_S\inv - C_0 - \ol C_0$, and therefore 
$$
\deg (N_{S_1/Z_1}|_{C_i}) = K_S\inv\cdot C_i - C_0\cdot C_i
- \ol C_0\cdot C_i  = 0 - 1-0 = -1.
$$
Hence from the last exact sequence, we obtain $N_{C_i/Z_1}
\simeq \ms O(-2)\oplus\ms O(-1)$.
From this, we obtain 
$$N_{X_1/Z_1}|_{C_i}\simeq X_1|_{C_i}\simeq
(X_1+\ol X_1)|_{C_i}\simeq
2S_1|_{C_i}\simeq 2(S_1|_{S_1})|_{C_i}
\simeq 2(K_S\inv - C_0 - \ol C_0)|_{C_i}\simeq
\ms O_{C_i}(-2).$$
Hence by using the inclusion 
$C_i\subset X_1\subset Z_1$,
 we can derive
$N_{C_i/X_1}\simeq\ms O(-1)$.
Thus each $C_i$ is a $(-1)$-curve on $S_1$.
Hence we can contract these four curves on $S_1$,
and consequently the self-intersection number
of the curve $E_0|_{X_1}$  becomes
$2(2-n) + 4 = 2(4-n)$ on the blowing-down.

By Lemma \ref{lemma:fl}, 
as $X\simeq F-\aaa_{5,6,\dots,n}$, we have $X.\,l_i=2$,
and the intersection 
$X\cap l_i$ consists of two points for any 
$i\in\{5,6,\dots,n\}$ in a generic situation, 
and each of the two intersection points is joined with
the elliptic curve $E_0|_{X_1}$ by a $(-1)$-curve.
Moreover if we blow-down all these $2(n-4)$-curves,
the self-intersection number of the curve
$E_0|_{X_1}$ becomes precisely zero.
If $X\cap l_i$ happened to consist of one point (for some $i$),
this point is an ordinary double point, and by resolving it
we obtain a $(-2)$-curve.
This curve has exactly two fixed points of the $\CC^*$-action,
and one of the two points intersects $(-1)$-curve which joins
the $(-2)$-curve with the elliptic curve $E_0|_{X_1}$.
Therefore again by contracting these
$(-1)$-curve and $(-2)$-curve subsequently, 
we again obtain that the contribution from the fixed twistor line $l_i$
for increasing the self-intersection number
of the elliptic curve $E_0|_{X_1}$ is 
again 2.

Thus all the intersections $X\cap l_i$ are joined by 
exceptional curves for obtaining the minimal model $X_{\rm{min}}$.
On the other hand, since a general member of the pencil $|2F|$ is
irreducible, the special member $X+\ol X$ is  connected,
and therefore $X\cap \ol X$ has to be a curve, which is of course
$\CC^*$-invariant.
The classification result \cite[Klassifikationssatz]{Ha95} means that this curve is either a smooth elliptic curve
or a cycle of rational curve,
depending on whether $X_1$ is birational to a Hopf surface
or a parabolic Inoue surface respectively.
But the latter cannot happen since $X\cap \ol X$ cannot have
any fixed point of the $\CC^*$-action anymore.
Hence $X_1$, and so $X$, is birational to a Hopf surface, 
as asserted.\proofend

\medskip

The intersection $X\cap \ol X$ has the following characteristic
properties:

\begin{proposition}
The intersection  $X\cap \ol X$ is a smooth elliptic curve, which
is a single orbit of the $\CC^*$-action on $Z$.
Moreover, $X\cap \ol X$ is homologous to zero.
\end{proposition}

\proof
From the above proof of Theorem \ref{thm:1},
we just need to show that 
the curve $X\cap \ol X$ is homologous to zero.
For this,
it suffices to show 
\begin{align}\label{hom0}
X.\ol X. y = 0\quad{\text{for any $y\in H^2(Z,\ZZ)$}}.
\end{align}
Let $\xi_1,\xi_2,\xi_3,\xi_4$ be elements
of $H^2(n\CP^2,\ZZ)$ such that 
$\{\xi_i\set 1\le i\le n\}$ is an orthonormal basis
of $H^2(n\CP^2,\ZZ)$, and put $\aaa_i:=\pi^*\xi_i$.
(We have already chosen $\xi_i$ for $5\le i \le n$ right before Lemma \ref{lemma:van1}.)
Then by \cite[p.\,31]{KK92} or \cite{P92}, 
$H^2(Z,\ZZ)$ is generated by the following $(n+1)$ elements:
$$
\frac 12 F + \frac12\sum_{1\le i\le n}\aaa_i
\quad{\text{and}}\quad \aaa_i,\,1\le i\le n.
$$
As $X= F - \aaa_{5,6,\dots,n}$ in $H^2(Z,\ZZ)$, we have
$$
X.\,\ol X = F^2 - ( \aaa_5^2 + \aaa_6^2 + \dots + \aaa_n^2)
\quad{\text{in}}\quad H^4(Z,\ZZ).
$$
Therefore if $y=\aaa_i$ for some $i$,
we have $X.\ol X.y = 0$ since $F^2.\aaa_i = 0$ and $\aaa_j^2.\aaa_i=0$ 
for any $i$ and $j$.
Moreover, we have
\begin{align*}
X.\ol X.\Big(
F + \sum_{1\le i\le n}\aaa_i
\Big)&= 
F^3 - ( \aaa_5^2 + \aaa_6^2 + \dots + \aaa_n^2).\,F\\
&= (8-2n) -(n-4)\cdot(-2)=0.
\end{align*}
Hence we obtain \eqref{hom0}.
\proofend


%


\subsection{Elimination of the base locus of the pencil $|2F|$}
As in Proposition \ref{prop:base35}, the anti-canonical system
$|2F|$ of the present twistor spaces has  base curves.
In this subsection, we give a sequence of blowups which completely
eliminates the base locus.
This gives an explicit holomorphic model for the algebraic reduction of $Z$. 
We begin with the following easy

\begin{proposition}\label{prop:bs2}
Denoting $C_1,\dots, C_4$ for the strict transforms of the 
original curves in $Z$ into the blowup $Z_1$,
for the line bundle $L_1= \mu_1^*(2F) - 2E_0 - 2\ol E_0$,
 we have
\begin{align}\label{bs3}
\Bs\,|L_1| =
(C_1\cup C_2\cup C_3\cup C_4) 
\sqcup
(\ol C_1\cup \ol C_2\cup\ol C_3\cup \ol C_4) .
\end{align}
Moreover, a general member of the pencil $|L_1|$
is  non-singular.
\end{proposition}

\proof
As  in \eqref{bs2}, if $X_1$ is the strict transform of
$X$ into $Z_1$ as before, we have 
$X_1|_S = C_{1,2,3,4}$.
Moreover, the pencil $|L_1|$ is generated by $2S_1$ and 
$X_1+ \ol X_1$.
These imply \eqref{bs3}.
Therefore by Bertini's theorem, a general member of the pencil 
$|L_1|$
can have singularity only on $C_{1,2,3,4}$
and $\ol C_{1,2,3,4}$.
Since $C_i$-s are disjoint, $X_1|_S = C_{1,2,3,4}$ already means that 
$X_1$ is smooth on $C_{1,2,3,4}$.
Similarly, $\ol X_1$ is smooth on $\ol C_{1,2,3,4}$.
These mean that a general member of $|L_1|$ is also
smooth on $C_{1,2,3,4}$ and $\ol C_{1,2,3,4}$.
\proofend

\medskip
Thus the double curves of a general member of the pencil $|2F|$
are completely resolved by the blowup $\mu_1$.
For the resolved divisors, we have the following

\begin{proposition}
Let $Y\in |2F|$ be a general anti-canonical divisor of $Z$,
and $Y_1$ its strict transform into $Z_1$.
Then we have
$$
K_{Y_1}\simeq -(E_0+\ol E_0)|_{Y_1},\quad
K_{Y_1}^2 = 4(2-n).
$$
Moreover, $Y_1$ contains the curves $C_i$ and $\ol C_i$
($i\in\{1,2,3,4\}$) as  $(-1)$-curves.
\end{proposition}

Since this can be shown easily by using adjunction formula
and the standard exact sequence associated to the inclusions
$C_i\subset S_1\subset Z_1$ and $C_i\subset Y_1\subset Z_1$,
we omit a proof.

%

Next we investigate the way how general members
of the pencil $|2F|$ intersect along the curve $C_i$
for $i\in\{1,2,3,4\}$:

\begin{proposition}\label{prop:tan}
A general member of the  pencil $|2F|$
is tangent to the divisor $X$ along
the curve $C_i\cup\ol C_i$ for any $i\in\{1,2,3,4\}$
to the second order.
\end{proposition}

\proof
Let $x$ be a defining section of the divisor $X\in
|F-\aaa_{5,6,\dots,n}|$.
Similarly, let $s\in H^0(F)$ be a non-zero element,
so that $(s) = S$.
Then since $S$ is non-singular, and
since $X$ is non-singular at least on 
 the curve $C_i$ ($i\in\{1,2,3,4\}$)
except possibly at the intersection point $C_i\cap C_0$,
from \eqref{rest1},
we can take $s$ and $x$ as a part of local coordinates
of $Z$ around the point.
Then since the pencil $|2F|$ is generated by 
$2S$ and $X+\ol X$, 
a defining section of a general member of the pencil
is of the form 
\begin{align}\label{tangent2}
x\ol x = cs^2,\quad
c\in\CC^*.
\end{align}
Moreover, from \eqref{rest1},
the divisor $\ol X$ does not intersect
$C_{1,2,3,4}$, and therefore 
$\ol x\neq 0$ for any point of $C_{1,2,3,4}$.
From the equation \eqref{tangent2},
this  means that a general member of $|2F|$ is tangent
to $X$ along $C_{1,2,3,4}$ to the second order.
\proofend

%

\medskip
Note however that the intersection of the divisors $S$ and $X$ 
is transverse along $C_{1,2,3,4}$.

Proposition \ref{prop:tan} means that,
 even after we blow up $Z_1$ at the base curves \eqref{bs3},
the strict transform of the pencil still has base curves
(which is identified with the original base curve \eqref{bs3}).
Because the tangential order of general members 
of the pencil is two,
it is  not difficult to show that 
if we blow up the  base curves \eqref{bs3},
then the resulting pencil has no base point.
In summary, we need the next 3 steps for the elimination of the base locus
of $|2F|$ on $Z$:
\begin{enumerate}
\item[1.]
The blowup $\mu_1:Z_1\to Z$ at the curves $C_0\cup\ol C_0$.
This resolves the  double curves  of  general members of 
the pencil $|2F|$.
\item[2.]
The blowup $Z_2\to Z_1$ at the curves $C_1\cup C_2\cup C_3\cup
C_4$ and the conjugate curves.
This transforms the touching situation of general members
along these eight curves into transversal intersections.
\item[3.]
The blowup $Z_3\to Z_2$ at the transversal intersections of 
general members of the pencil.
This blowup  separates the intersections completely.
\end{enumerate}

\section{Existence of a twistor space on $5\CP^2$ with a pencil of K3 
surfaces}
In this section, we show that, on $5\CP^2$, there exists a twistor space 
$Z$ satisfying $a(Z)=1$ whose general fiber of the algebraic reduction is
a K3 surface.
A proof proceeds in a similar way to the case of elliptic ruled surfaces
employed in Sections 2 and 3,
but the proof is pretty easier in that we do not need  complicated
calculations for  cohomology groups like we did in Section 3.1. 

\begin{figure}
\includegraphics{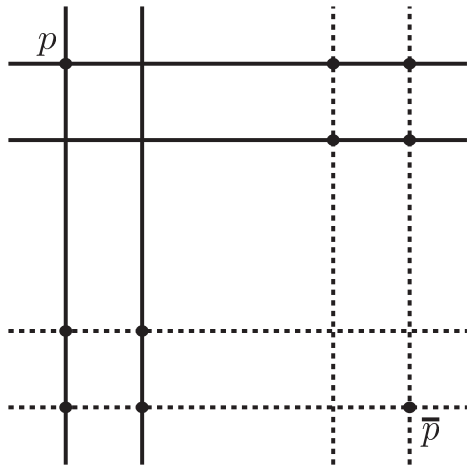}
\caption{Blown up points on $\CP^1\times\CP^1$ to construct the 
surface $S$}
\label{fig:K3}
\end{figure}

First we construct a rational surface $S$ with $c_1^2(S)=-2$ 
which will be included in a twistor space on $5\CP^2$ as a real member of the linear
system $|F|$.
For this, we start from the product surface $\CP^1\times \CP^1$.
On it we put a real structure that is the product of the complex conjugation
and the anti-podal map.
We mean by a $(k,l)$-curve on $\CP^1\times\CP^1$ a curve of bidegree $(k,l)$.
We first take distinct two non-real $(1,0)$-curves.
We suppose that these are {\em not} mutually conjugate to each other.
Next we take distinct two $(0,1)$-curves, which are {\em not} mutually conjugate to each other.
(Each one cannot be real from the choice of the real structure
on $\CP^1\times\CP^1$.)
Let $D$ be the sum of all these four curves.
This is a reduced, non-real $(2,2)$-curve.
In Figure \ref{fig:K3}, the components of $D$ are displayed as
solid lines.
Then the conjugate curve $\ol D$ is also a reduced, non-real $(2,2)$-curve,
whose components are displayed as dashed lines in Figure \ref{fig:K3}.
The intersection $D\cap \ol D$ consists of 8 points.
Next we choose any one of the four double points,
say $p$, of the curve $D$.
Let 
$$
\eee:S\to \CP^1\times\CP^1
$$
be the blowup at the 10 points $(D\cap \ol D)\cup\{p,\ol p\}$.
Let $\sigma$ be the real structure on $S$ which is 
a natural lift of that on $\CP^1\times\CP^1$.
For this surface $S$, we have the following

\begin{proposition}\label{prop:S5}
Let $S$ be the rational surface with $c_1^2(S)=-2$ constructed as above. Then $H^0(mK_S\inv)=0$ if $m$ is odd and 
$H^0(mK_S\inv)\simeq \CC$ if $m$ is even.
In particular, $\kappa\inv(S)=0$ for the anti-canonical dimension
of $S$.
\end{proposition}

\proof
Let $C$ be the strict transform of the curve $D$ into $S$,
and $C_5:=\eee\inv(p)$ the exceptional curve over $p$.
Then a double point $p$ of $D$ is 
separated by $\eee$, and 
$C$ is a chain of  four rational curves.
The self-intersection numbers 
of these four rational curves are readily seen to be
\begin{align}\label{chain1}
-3,-2,-2,-3
\end{align}
respectively.
Of course the same is true for the other chain $\ol C$.
Note that $C\cap \ol C=\emptyset$.
Moreover, noting $C\simeq \eee^*D - 2C_5$, we easily obtain linear equivalences
\begin{align}\label{chain2}
C\simeq K_S\inv - C_5 + \ol C_5
\quad{\text{and}}\quad
\ol C\simeq K_S\inv + C_5 - \ol C_5.
\end{align}
From this we obtain $2K_S\inv\simeq C + \ol C$.
Thus the sum $C+\ol C$ belongs to $|2K_S\inv|$.
Since the intersection matrices for $C$ and $\ol C$ may be
readily seen to be negative definite from \eqref{chain1},
by a property of the Zariski decomposition,
we obtain that $H^0(2mK_S\inv)\simeq\CC$ for 
any $m>0$.
The other assertion $H^0((2m-1)K_S\inv)=0$ can be shown in 
a similar way to the same assertion in Proposition \ref{prop:rs2} (ii).
So we omit the detail.
\proofend

\medskip
Similarly to what we did in Section 3.1, let $\aaa_5\in H^2(Z,\ZZ)$ be
the cohomology class determined from the curve $C_5$ on $S$.
So we have $\aaa_5|_S = C_5-\ol C_5$ in $H^2(S,\ZZ)$.
The next proposition amounts to Propositions \ref{prop:Z1} and \ref{prop:nonvan1}
in the case for twistor spaces studied in Section 3, but 
the proof is much easier basically because $n=5$.
\begin{proposition}\label{prop:Z16}
Let $Z$ be a twistor space on $5\CP^2$ which includes
the surface $S$ as a real member of $|F|$.
Then we have:
\begin{itemize}
\item[\em (i)]
The system $|F|$ consists of the single member $S$.
\item[\em (ii)]
The system $|F-\aaa_5|$ consists
of a single member.
If $X$ denotes the member, $X$ is irreducible
and $X|_S = C$,
where $C$ is the above chain of four curves on $S$.
\item[\em (iii)]
The anti-canonical system $|2F|$ of $Z$ is a pencil,
and is generated by a non-reduced member $2S$ and 
the reducible member $X + \ol X$.
\item[\em (iv)]
$a(Z)=1$, and the algebraic reduction of $Z$ is induced by the 
meromorphic map associated to
the pencil $|2F|$.
\end{itemize}
\end{proposition}

\proof
The item (i) immediately follows from $|K_S\inv|=\emptyset$
in Proposition \ref{prop:S5}.
For (ii), we consider the standard exact sequence
\begin{align}\label{les79}
0 \lras \ms O_Z(-\aaa_5)
\lras F-\aaa_5
\lras (F -\aaa_5)|_S\lras 0.
\end{align}
Then from the fact that $(C_5)^2_S=-1$ on $S$,
exactly in the same way to the final part of the proof of 
\cite[Proposition 3.3]{HonCrelle2},
we obtain, by Riemann-Roch formula and Hitchin's vanishing
theorem, that 
$$
H^1(\ms O_Z(-\aaa_5)) = 0.
$$
Hence from the cohomology exact sequence of \eqref{les79}
we obtain an isomorphism
$H^0(F-\aaa_5) \simeq H^0((F-\aaa_5)|_S)$.
Moreover, from an isomorphism in \eqref{chain2}, we have
\begin{align*}
(F-\aaa_5)|_S & \simeq  K_S\inv - (C_5 - \ol C_5) \simeq C.
\end{align*}
Hence $H^0((F-\aaa_5)|_S)\simeq H^0(\ms O_S(C))$.
Irreducibility of the unique member $X$ of $|F-\aaa_5|$
again follows from $|F|=\{S\}$ as in the proof of Proposition 
\ref{prop:X1} (ii).
This proves the assertion (ii) of the proposition.

For (iii) and (iv), since $2S$ and $X + \ol X$ are 
linearly independent
anti-canonical divisors of $Z$, we have 
$h^0(2F)\ge 2$.
Further from the exact sequence
$0 \to F \to 2F \to 2K_S\inv\to 0$ and
recalling $h^0(F) = 1$ and $h^0(2K_S\inv) = 1$
from Proposition \ref{prop:S5}, we obtain
$h^0(2F)\le 2$. So $h^0(2F) = 2$.
Hence $a(Z)\ge 1$.
On the other hand, since $\kappa\inv(S)=0$, we have
$a(Z) \le 1 + \kappa\inv(S)=1$. So $a(Z)=1$.
Finally the second assertion in (iv) can be obtained in the same way
to the final part of the proof of Proposition \ref{prop:Z1}
(which is at the end of Section 3.1).
\proofend

%
%
\medskip
Next we will determine structure of a general member
of the pencil $|2F|$.
For this we first see

\begin{proposition}\label{prop:Bs16}
Let $Z$ be a twistor space on $5\CP^2$ which includes
the surface $S$ as a real member of $|F|$
as in Proposition \ref{prop:Z16}.
Then we have $\Bs\,|2F|= C\cup \ol C$.
\end{proposition}

\proof
Since the pencil $|2F|$ is generated by 
the two divisors $2S$ and $X+\ol X$ as
in Proposition \ref{prop:Z16}, 
we have $\Bs\,|2F|= S\cap (X\cup\ol X)$.
This is the chains $C\cup\ol C$ by Proposition \ref{prop:Z16} (ii).
\proofend

\medskip
Thus the chains $C\cup\ol C$ play a similar role
to the curves \eqref{isom0} in the case 
of the twistor spaces discussed in Section 3.
Main difference is that in the present case, the restriction
$X|_S$ is a reduced divisor on $S$, while
it is non-reduced
for the other type of twistor spaces as in \eqref{rest1}.
This yields the following difference on the structure 
of a general member of the pencil $|2F|$:

\begin{theorem}\label{thm:2}
Let $Z$ be a twistor space on $5\CP^2$ which includes
the surface $S$ as a real member of $|F|$
as in Proposition \ref{prop:Z16}.
Then a general member of the pencil $|2F|$ is birational
to a K3 surface.
In particular, a general fiber of the algebraic reduction
of $Z$ is birational to a K3 surface.

\end{theorem}

\proof
Recalling that $\Bs\,|2F| = C\cup \ol C$, let 
$$
\mu_1:Z_1\to Z
$$
be the blowup of $Z$ at the two chains $C\cup \ol C$,
and $E_1\cup\ol E_1$ the exceptional divisor over $C\cup\ol C$.
The restrictions $\mu_1|_{E_1}:E_1\to C$ and $\mu_1|_{\ol E_1}:\ol E_1\to \ol C$ are $\CP^1$-bundle maps, and $Z_1$ has one ordinary
double point
 on each of the fibers over a double point of $C$ and $\ol C$.
Then as $X|_S=C$ as in Proposition \ref{prop:Z16},
the divisors $X$ and $S$ are already completely separated in $Z_1$.
Define
$$
L_1:=\mu_1^*(2F) - E_1 - \ol E_1.
$$
Then from Proposition \ref{prop:Bs16} we have an isomorphism $|L_1|\simeq 
|2F|$.
Since the pencil $|2F|$ is generated by $2S$ and $X+\ol X$,
the pencil $|L_1|$ is generated by, this time, 
$$
2S_1+E_1+\ol E_1\quad
{\text{and}}\quad X_1+\ol X_1,
$$
where $S_1, X_1$ and $\ol X_1$ are strict transforms of $S, X$ and $\ol X$
into $Z_1$ respectively.
Since $X$ and $S$ are already separated in $Z_1$ as above,
 we have 
\begin{align}\label{bs62}
\Bs\,|L_1| = (X_1\cap E_1)\cup (\ol X_1\cap \ol E_1).
\end{align}
Both $X_1\cap E_1$ and $\ol X_1\cap \ol E_1$ are chains of 
curves which are isomorphic to the original chains $C$ and $\ol C$
respectively by $\mu_1$ .

As above, $Z_1$ has one ordinary double point over
each of the double points of the chains $C$ and $\ol C$.
Let $q\in C$ be any one of the three double points of the chain $C$,
and suppose that the divisor $X$ (in $Z$) is singular at $q$.
Then the strict transform 
$X_1$ passes the singular point of $Z_1$ over the point $q$ 
iff $X$ has
a singularity at $q$, as one can readily be seen from computations using local coordinates.
Regardless of whether $q\in \Sing X$ or not, let 
$$
\mu_2:Z_2\to Z_1
$$
be the blowup at the base curve \eqref{bs62} of the pencil $|L_1|$.
Let $E_2$ and $\ol E_2$ be the exceptional divisors
over the chains $X_1\cap E_1$ and $\ol X_1\cap \ol E_1$
respectively.
While $Z_2$ has singularities, the blowup $\mu_2$ already
eliminates the base locus of the pencil $|2F|$ and $|L_1|$.
Namely, if we put
$$
L_2:=\mu_2^*L_1 - E_2 - \ol E_2,
$$
then we have $\Bs\,|L_2| = \emptyset$.

Type of the singularities of $Z_2$ depends
on whether $q\in \Sing X$.
If $q\not\in \Sing X$, 
over the point $q$, the variety $Z_2$ has exactly two 
ordinary double points, one of which appears from
the first blowup $\mu_1:Z_1\to Z$,
and the other of which appears from the second blowup
$\mu_2:Z_2\to Z_1$.
So both of the singularities admit a small resolution.
On the other hand,
if $q\in \Sing X$, the ordinary double point of 
$Z_1$ splits by the blowup $\mu_2:Z_2\to
Z_1$ into two singularities,
both of which can be locally written by an equation
$xy = zw^2$ in $\CC^4(x,y,z,w)$ in 
local coordinates.
This singularity also admits a small resolution, which may be concretely 
obtained by blowing-up the plane $\{x=w=0\}$
in the last expression.
Taking this small resolution for each of the two
singularities,
both  are transformed into 
an ordinary double point.
So these again admit a small resolution.
Let 
$$
\mu_3:Z_3\to Z_2
$$
be any one of the small resolutions obtained this way,
which preserves the real structure.
We put
$$
L_3:=\mu_2^*L_2.
$$
Then since $\Bs\,|L_2| = \emptyset$ as above,
we have $\Bs\,|L_3| = \emptyset$.
Let 
$$
\Phi_3:Z_3\lra \CP^1
$$
be the morphism induced by the pencil $\Bs\,|L_3|$.
From the construction, there is a natural identification
between the two pencils $|2F|$ on $Z$ and $|L_3|$ on $Z_3$,
and members of $|2F|$ are naturally identified with
fibers of the morphism $\Phi_3$.
Hence in order to show the theorem,
it is enough to see that a general fiber of $\Phi_3$ 
is a K3 surface.

For this, we write $\mu$ for the composition $\mu_1\circ\mu_2\circ\mu_3$,
and in the following $E_1$ and $E_2$ mean the strict transforms of 
the original exceptional divisors into $Z_3$.
Then from the usual transformation formula for the canonical bundle
under a blowup,
outside the 6 singularities of the chains $C\cup\ol C$ in $Z$,
we have
\begin{align}\label{ant}
K_{Z_3}\inv\simeq \mu^* (K_Z\inv) - E_1 - \ol E_1
-2 E_2 - 2\ol E_2.
\end{align}
Since this is valid outside codimension two locus, 
by Hartogs' theorem, \eqref{ant} is valid on the whole of $Z_3$.
But the right-hand side of \eqref{ant} is exactly the line bundle
$L_3$.
Hence a general fiber of $\Phi_3$ is an anti-canonical divisor on 
$Z_3$.
Hence for a general fiber $Y_3$ of $\Phi_3$, we have 
$K_{Y_3}\simeq \ms O$.
To complete a proof, it suffices to show $H^1(Y_3,\ms O)=0$.
For this we consider the standard exact sequence
\begin{align*}
0\lras \ms O_{Z_3}(-Y_3)
\lras \ms O_{Z_3} \lras \ms O_{Y_3}\lras 0.
\end{align*}
As $Y_3\in |K_{Z_3}\inv|$ as above, 
we have $\ms O_{Z_3}(-Y_3)\simeq K_{Z_3}$ for the 
first term.
Hence $H^1(\ms O_{Z_3}(-Y_3))$ is dual to 
$H^2(\ms O_{Z_3})$.
The latter is of course isomorphic to $H^2(\ms O_Z)$ from 
birational invariance.
Moreover, we have $H^2(\ms O_Z)=0$ since this is again dual 
to $H^1(K_Z)$, which vanishes by Hitchin's vanishing theorem.
Thus $Y_3$ is a K3 surface, as asserted.
\proofend

\medskip
Of course, the morphism $\Phi_3:Z_3\to\CP^1$ obtained above may be regarded
as a holomorphic family of K3 surfaces.
The two fibers corresponding to the divisors $2S$ and $X+\ol X$ are
degenerate fibers of this fibration.
Thus the twistor spaces in Theorem \ref{thm:2} provide
degeneration of K3 surfaces under a non-K\"ahler setting.
In this regard, it might be interesting to see if how one can 
make these degenerations to be semi-stable,
and after that, which type of degenerations occur \cite{FM83,Ni88}.
We note that the morphism $\Phi_3$ is induced by the anti-canonical 
system of $Z_3$, and therefore the canonical bundle
of $Z_3$ is trivial 
around the singular fibers.

\section{Existence of the twistor spaces and dimension of the moduli spaces}

\subsection{Existence}
In this subsection we show that the twistor spaces 
investigated in Sections 3 and 4 actually exist.
Recall that both types of the twistor spaces are characterized by 
the complex structure of the (unique) member $S$ of $|F|$.
Thus it is enough to see that for such a
prescribed surface $S$,
there exists a twistor space $Z$ having $S$
as a member of $|F|$.
We show this by deforming known example of a twistor space
$Z_0$ on $n\CP^2$ having a divisor belonging to $|F_{Z_0}|$,
in such a way that the divisor is preserved and deformed into 
the prescribed surface $S$.

Let $n\ge 5$ be any integer. 
First we show the existence of a twistor space $Z$ on $n\CP^2$
enjoying the assumptions in Theorem \ref{thm:1}.
For this, we consider the twistor space of a self-dual metric
on $n\CP^2$ constructed by Joyce \cite{J95}.
These metrics admit a $T^2$-action,
and its natural lift to the twistor space generates
a holomorphic $(\CC^*\times\CC^*)$-action.
By the result of Fujiki \cite{F00}, the closure of a generic orbit
of the $(\CC^*\times\CC^*)$-action is a non-singular toric surface,
and the toric surface is uniquely determined
from (the isomorphism class of) the $T^2$-action on $n\CP^2$.
Conversely, the the toric surface uniquely determines
the $T^2$-action on $n\CP^2$ of a Joyce metric.
Thus in order to specify the $T^2$-action of the Joyce metric,
it is enough to specify the toric surface.

\begin{figure}
\includegraphics{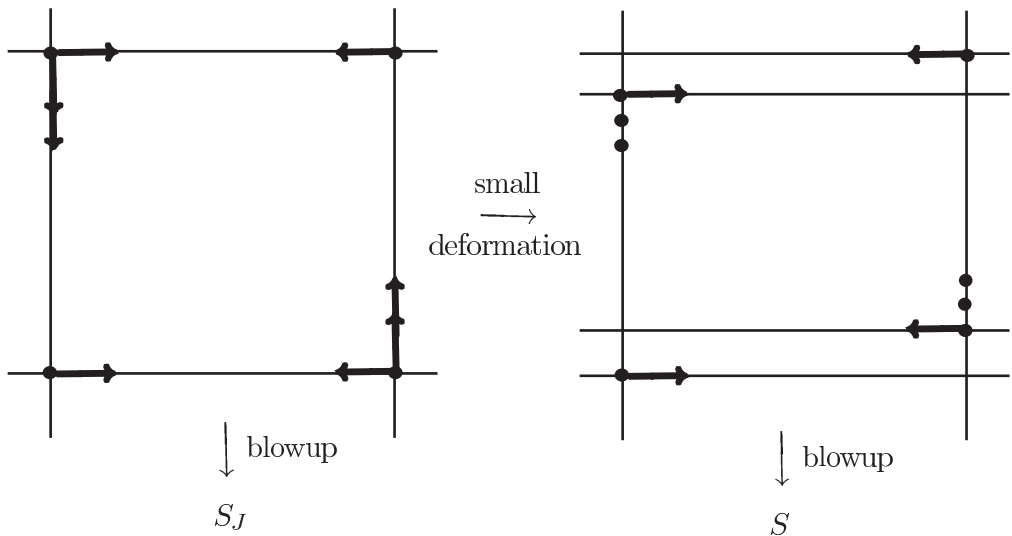}
\caption{deformation from $S_J$ to $S$ for the case of elliptic ruled surface}
\label{fig:dfm}
\end{figure}

For specifying the toric surface,
we consider the product
$\CP^1\times\CP^1$ equipped with the real structure 
used in the beginning of Section 4.
We again take a reducible $(2,2)$-curve consisting two 
$(1,0)$-curves and two $(0,1)$-curves, but this time
we take them in a way that the $(2,2)$-curve becomes real.
We make an iterated blowup at each of the four double points 
in the way as displayed in Figure \ref{fig:dfm}.
Here, the arrows mean to blowup the exceptional curve
at the point which corresponds to the direction indicated
by the arrow.
Let $S_J$ be the toric surface obtained by this blowup.
Since we blowup $\CP^1\times\CP^1$ precisely $2n$ times,
we have $K^2=8-2n$ for $S_J$.
Also $S_J$ admits a natural real structure as a lift from
$\CP^1\times\CP^1$.
Then there exists a Joyce metric on $n\CP^2$ whose twistor space
has the toric surface $S_J$ as a $T^2$-invariant member
of the system $|F|$.

We consider a deformation of $S_J$ preserving the real structure
in a way displayed in Figure \ref{fig:dfm}.
Namely we shift the pair of a point and one arrow on 
lower left
in the vertical direction, and at the same time, 
displace all vertical $(n-4)$ arrows on upper left in
the vertical direction to 
split them into distinct $(n-4)$ points.
We move the points on the right side as determined by
the real structure from the left side.

Let $S$ be the surface obtained from $\CP^1\times\CP^1$ 
by the above blowing up.
Note that if $n=4$, $S$ is nothing but the 
relatively minimal elliptic surface $S_0$  constructed 
in Section 2.
We also note that this deformation (of $S_J$ into $S$)
preserves a $\CC^*$-action; we just need to consider 
the product of a standard $\CC^*$-action on 
horizontal direction and the trivial $\CC^*$-action on
the vertical direction.
In summary the surface $S$ in Theorem \ref{thm:1}
is obtained from the toric surface $S_J$ by a 
small deformation preserving the real structure and a $\CC^*$-action.

Once this is obtained, it is almost automatic to show the 
existence of a twistor space $Z$ on $n\CP^2$ having $S$ as a
real member of the system $|F|$.
For this, letting $Z_0$ be 
the twistor space of a Joyce metric having 
the divisor $S_J$ as above,  we consider deformation of the pair $(Z_0,S_J)$ preserving the $\CC^*$-action and the real structure.
Then  because $ H^2(Z_0,\Theta_{Z_0}(-S_J)) = 0$ from \cite[Lemma 1.9]{C91-2},
the divisor $S_J$ is co-stable in $Z_J$ \cite[Theorem 8.3]{Hor76}.
This proves the existence of a twistor space on $n\CP^2$
satisfying the assumptions in Theorem \ref{thm:1}.

\begin{remark}{\em
Because the twistor spaces of Joyce metrics do not contain
a non-algebraic surface,
one may wonder from where a VII surface in the present twistor space comes from.
An answer is that on the twistor space of a Joyce metric,
there exist two degree-one divisors $D$ and $D'$ 
such that the sum $D+ D'$ belongs to the same cohomology class
as the divisor $X$,
and under the small deformation the sum becomes the (irreducible) VII surface $X$.
}
\end{remark}

\begin{figure}
\includegraphics{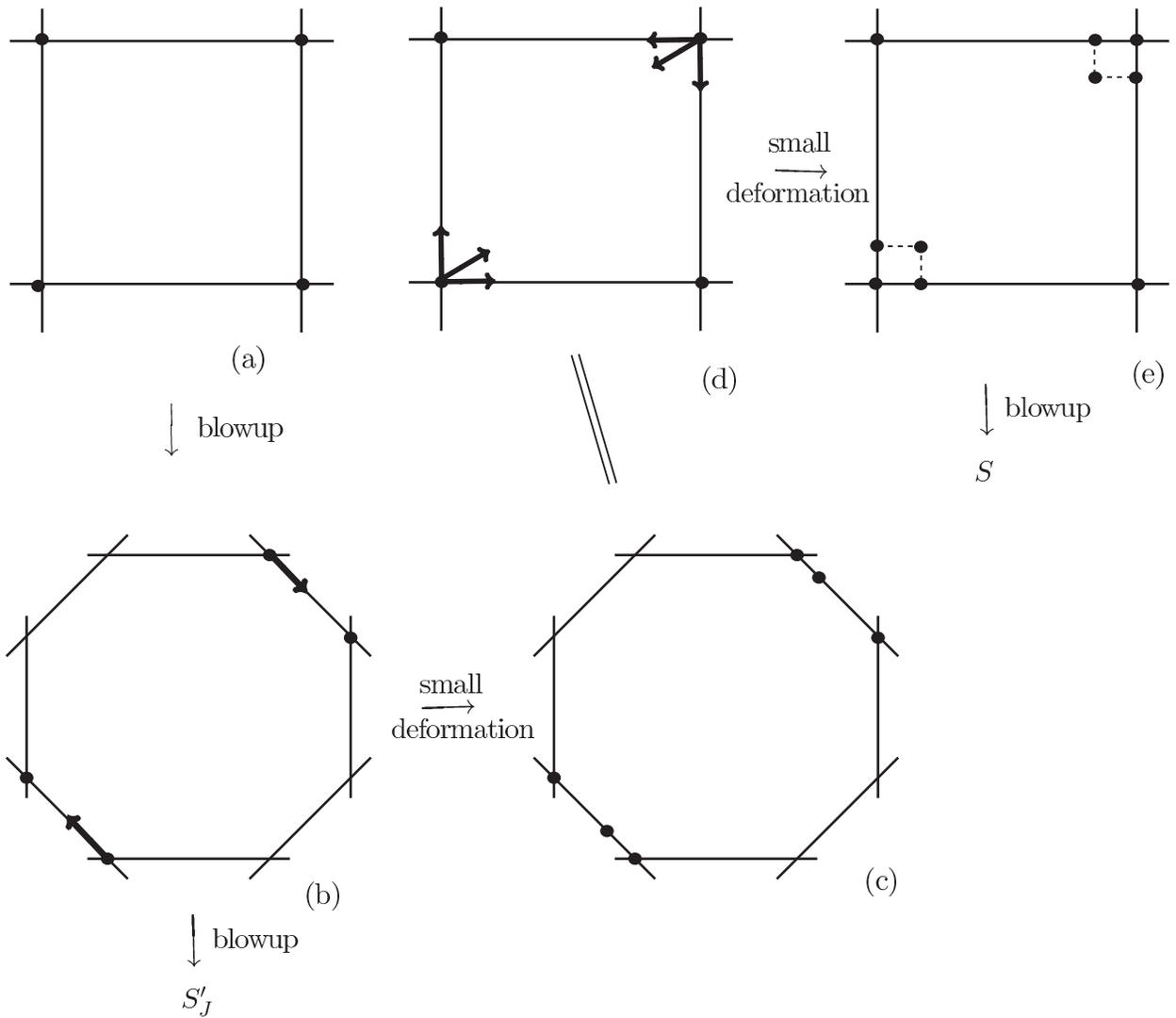}
\caption{deformation from $S'_J$ to $S$ for the case of K3 surface}
\label{fig:dfm2}
\end{figure}

The existence of a twistor space satisfying the assumptions
in Theorem \ref{thm:2} can be shown in a similar way,
with slightly more effort.
We again start from a toric surface,
for which we denote $S'_J$ this time, which is contained 
in the twistor space of a Joyce metric as a $T^2$-invariant
member of $|F|$.
This toric surface $S'_J$ is obtained from $\CP^1\times\CP^1$ 
by blowing up 10 times as displayed in Figure \ref{fig:dfm2} (a) and (b).
Namely we first blowup $\CP^1\times\CP^1$ at the 4 
double points (Figure \ref{fig:dfm2} (a)) of the reducible $(2,2)$-curve,
and then blowup the resulting surface 6 times 
at the points and arrows in Figure \ref{fig:dfm2} (b).
Let $S'_J$ be the resulting toric surface.
We have $K^2 = -2$ for this toric surface.
Next we move each arrow to a point in that direction (
as indicated in Figure \ref{fig:dfm2} (c)).
By blowing up these 6 points we obtain a non-toric surface as a small deformation
of $S'_J$.
This surface still admits a $\CC^*$-action.
In terms of the initial surface $\CP^1\times\CP^1$,
the way of blowing up for obtaining this surface with $\CC^*$-action can be
indicated as in Figure \ref{fig:dfm2} (d).
Next we deform this configuration in a way 
that the three arrows from the point on lower left
move along  lines in each direction, in a way
as displayed in Figure \ref{fig:dfm2} (e).
Here, `line' means a line in the usual sense,
viewing $\CP^1\times\CP^1$ as a compactification
of $\CC\times\CC=\CC^2$.
Let $S$ be the surface obtained from $\CP^1\times\CP^1$
by blowing up the  10 points in Figure \ref{fig:dfm2} (e).
Then this is exactly the surface $S$ constructed at the beginning of
Section 4.
Hence we have seen that the surface $S$ is obtained from the toric surface
$S'_J$ by first deforming it preserving a $\CC^*$-action,
and then further deforming the resulting surface by a small deformation.
Since the cohomology vanishing property is preserved
under small deformation, by the same reasoning to the above case,
we again obtain that 
there actually exists a twistor space $Z$ on $5\CP^2$ 
having the surface $S$ as a member of $|F|$.
Thus we have shown the existence
of a twistor space satisfying the assumption of Theorem \ref{thm:2}.

\subsection{Dimension of the moduli spaces}
In this subsection we calculate dimension of the moduli spaces
of the twistor spaces discussed in Sections 3 and 4.
First let $n\ge 5$ and  $Z$ be the twistor space on $n\CP^2$ studied in Section 3, and $S\in|F|$ the divisor which was supposed to
exist.
As in Proposition \ref{prop:Z1}, the divisor $S$ is the unique member
of the system $|F|$.
Hence if $(Z',S')$ is  another pair of a twistor space and 
the divisor on it enjoying the assumption in Theorem 
\ref{thm:1}, and if $S\not\simeq S'$, we have
$Z\not\simeq Z'$.
Moreover, as long as the twistor space $Z$ is obtained as
a small deformation of the twistor space of a Joyce metric
as in the last subsection,
we have $H^2(\Theta_Z(-S))=0$ by upper semi-continuity,
and therefore by \cite[Theorem 8.3]{Hor76}, any small deformation of $S$
may be realized by a deformation of the pair $(Z,S)$.
Hence in order to compute the dimension of the moduli space,
it is enough to compute the dimension of the moduli space
for which the complex structure of the divisor $S$ is {\em fixed.}

For this, we consider an exact sequence
\begin{align}\label{ses:87}
0\lras \Theta_{Z}(-S)\lras \Theta_{Z,S}\lras 
\Theta_S\lras 0.
\end{align}
As we have been assuming that the $\CC^*$-action on $S$ extends
to the twistor space $Z$, together with $H^2(\Theta_Z(-S))=0$ as above,
we get an exact sequence
\begin{align}\label{ses:88}
0\lras H^1(\Theta_{Z}(-S))\lras H^1(\Theta_{Z,S})\lras 
H^1(\Theta_S)\lras 0.
\end{align}
For computing the dimension $h^1(\Theta_{Z,S})$, we use the exact sequence
$
0\to \Theta_{Z,S}\to \Theta_Z\to K_S\inv\to 0.
$
As $H^0(K_S\inv) = 0$ by Proposition \ref{prop:rs2},
and $H^2(\Theta_{Z,S})=0$ from the upper semi-continuity,
we get an exact sequence
\begin{align}\label{ses75}
0\lras H^1(\Theta_{Z,S})\lras H^1(\Theta_Z)
\lras H^1(K_S\inv)\lras 0.
\end{align}
Then as $\chi(\Theta_Z) = 15-7n$, $h^0(\Theta_Z)=1$, and
$H^2(\Theta_Z)=0$ by upper semi-continuity,
we have
$
h^1(\Theta_Z) = 7n-14.
$
On the other hand, we readily have $h^1(K_S\inv) = 2n-9$.
Therefore from the  exact sequence \eqref{ses75} we 
obtain $h^1(\Theta_{Z,S}) = 5n-5$.
Moreover we have $\chi(\Theta_S) = 6-4n$ (as $S$ is a blow-up
of $\CP^1\times\CP^1$ at $2n$ points).
Hence since $h^0(\Theta_S)=1$ and $H^2(\Theta_S) = 0$,
we get $h^1(\Theta_S) = 4n-5$.
Therefore from the cohomology exact sequence of \eqref{ses:87},
under an assumption that the $\CC^*$-action on $S$ always
extends to that on the twistor space $Z$,
we obtain 
$$
h^1(\Theta_Z(-S)) = n.
$$
Therefore, 
we obtain that the moduli space of our twistor spaces with the complex structure
of the divisor $S$ fixed is $n$-dimensional.
On the other hand, as in Section 2, our surface $S$ is obtained from the relatively minimal
elliptic surface $S_0$ by blowing up $2(n-4)$ points on 
the rational curve $C_0$.
Moreover the moduli space for $S_0$ is identified with the moduli space of 
elliptic curves with a real structure, and therefore real 1-dimensional.
The contribution from the choice of the $2(n-4)$ points is 
real $2(n-4)$-dimensional,
so the moduli space for our surface $S$ is real $(2n-7)$-dimensional.
Thus under the above assumption on the extension of the $\CC^*$-action
from $S$ to $Z$, the moduli space of our twistor space is $(3n-7)$-dimensional.
Note that the assumption on the extension of $\CC^*$-action
might look somewhat strong, but as we argued in Section 3, 
we can derive the conclusion $a(Z)=1$ without using
the presence of $\CC^*$-action,
and just under the presence of the special divisor $S$.

We also remark that the natural $\CC^*$-action on the cohomology group $H^1(K_S\inv)\simeq
\CC^{2n-9}$ has exactly 1-dimensional subspace on which $\CC^*$-act trivially.
This seems to mean that the divisor $S$ disappears under a generic
small deformation of $Z$ which preserves $\CC^*$-action.
It might be interesting to ask what is the algebraic dimension of a twistor space obtained as 
these $\CC^*$-equivariant deformation. 

Next we compute the dimension of the moduli space of 
the twistor spaces investigated in Section 4.
First by the same reason to the case of twistor spaces in Section 3
discussed above,
the complex structure of the twistor space
$Z$ deforms if that of the divisor $S$ deforms.
But in the present case,
the complex structure of $S$ cannot be deformed,
 as is readily seen from the construction of the surface $S$.
 Further we have $H^0(K_S\inv)=0$,
 and so the exact sequence \eqref{ses75} is again valid,
 as long as the twistor space is obtained 
as a small deformation of the twistor space of a Joyce metric
as in  the last subsection.
Moreover, for these twistor spaces we have 
$h^1(\Theta_Z) = -\chi(\Theta_Z)=7\cdot 5 - 15 = 20$ and 
$h^1(K_S\inv) = 2\cdot 5 - 9 = 1$.
Hence from the exact sequence \eqref{ses75} we obtain 
$h^1(\Theta_{Z,S})= 19$.
Noting that the exact sequence \eqref{ses:88} is also valid
and noticing $h^1(\Theta_S) = -\chi(\Theta_S)=4\cdot 5 - 6 = 14$, we obtain
$h^1(\Theta_Z(-S))= 19-14=5$.
Thus the moduli space of the twistor spaces in Theorem \ref{thm:2} 
is 5-dimensional.

\end{document}